\documentclass[a4paper,12pt]{article}%
\usepackage{amsmath}
\usepackage{amssymb}
\usepackage{amsfonts}
\usepackage{amsthm,amscd}
\usepackage{indentfirst}
\usepackage{graphicx}%
\providecommand{\U}[1]{\protect\rule{.1in}{.1in}}
\def\figurename{Figure}
\makeatletter
\renewcommand{\fnum@figure}[1]{\figurename~\thefigure.}
\makeatother
\def\tablename{Table}
\makeatletter
\renewcommand{\fnum@table}[1]{\tablename~\thetable.}
\makeatother

\def \eop {\hbox{}\nobreak\hfill
\vrule width 2mm height 2mm depth 0mm
\par \goodbreak \smallskip}
\newtheorem{theorem}{Theorem}[section]
\newtheorem{lemma}[theorem]{Lemma}

\newtheorem{proposition}[theorem]{Proposition}
\newtheorem{definition}[theorem]{Definition}

\newtheorem{remark}[theorem]{Remark}
\numberwithin{equation}{section}

\begin{document}

\title{Time-consistent investment and consumption strategies under a general discount function}
\author{{I. Alia }\thanks{ Laboratory of Applied Mathematics, University Mohamed
Khider, Po. Box 145 Biskra (07000), Algeria. E-mail: ishak.alia@hotmail.com}
\and {F. Chighoub }\thanks{ Laboratory of Applied Mathematics, University Mohamed
Khider, Po. Box 145 Biskra (07000), Algeria. E-mail:
f.chighoub@univ-biskra.dz; chighoub\_farid@yahoo.fr}
\and {N. Khelfallah }\thanks{ Laboratory of Applied Mathematics, University Mohamed
Khider, Po. Box 145 Biskra (07000), Algeria. E-mail: nabilkhelfallah@yahoo.fr}
\and {J. Vives }\thanks{Departament de Matem\`{a}tiques i Inform\`{a}tica,
Universitat de Barcelona, Gran Via 585, 08007 Barcelona, Spain. E-mail:
josep.vives@ub.edu}}
\maketitle

\begin{abstract}
In the present paper we investigate the Merton portfolio management problem in
the context of non-exponential discounting, a context that give rise to
time-inconsistency of the decision maker. We consider equilibrium policies
within the class of open-loop controls, that are characterized, in our
context, by means of a variational method which leads to a stochastic system
that consists of a flow of forward-backward stochastic differential equations
and an equilibrium condition. An explicit representation of the equilibrium
policies is provided for the special cases of power, logarithmic and
exponential utility functions.

\end{abstract}

\textbf{Keys words}: Stochastic Optimization, Investment-Consumption Problem,
Merton Portfolio Problem, Non-Exponential Discounting, Time Inconsistency,
Equilibrium Strategies, Stochastic Maximum Principle.


\section{Introduction}

\bigskip

\hspace{0.4cm} \textit{Background}\bigskip

The common assumption in classical investment-consumption problems under
discounted utility is that the discount rate is assumed to be constant over
time which leads to the discount function be exponential. This assumption
provides the possibility to compare outcomes occurring at different times by
discounting future utility by some constant factor. But on the other hand,
results from experimental studies contradict this assumption indicating that
discount rates for the near future are much lower than discount rates for the
time further away in future. Ainslie, in \cite{Ains}, established experimental
studies on human and animal behaviour and found that discount functions are
almost hyperbolic, that is, they decrease like a negative power of time rather
than an exponential. Loewenstein \& Prelec in \cite{Loewenstein} show that
economic decision makers are impatient about choices in the short term but are
more patient when choosing between long-term alternatives, and therefore, a
hyperbolic type discount function would be more realistic.

Unfortunately, as soon as a discount function is non-exponential, discounted
utility models become time-inconsistent in the sense that they do not admit
the Bellman's optimality principle. Consequently, the classical dynamic
programming approach may not be applied to solve these problems. In light of
the non-applicability of dynamic programming approach directly, there are two
basic ways of handling time inconsistency in non exponential discounted
utility models. In the first one, under the notion of naive agents, every
decision is taken without taking into account that their preferences will
change in the near future. The agent at time $t\in\left[  0,T\right]  $ will
solve the problem as a standard optimal control problem with initial condition
$X(t)=x_{t}.$ If we suppose that the naive agent at time $0$ solves the
problem, his ot her solution corresponds to the so-called pre-commitment
solution, in the sense that it is optimal as long as the agent can pre-commit
his or her future behavior at time $t=0$. Kydland \& Prescott in
\cite{Kydland} indeed argue that a pre-committed strategy may be economically
meaningful in certain circumstances. The second approach consists in the
formulation of a time-inconsistent decision problem as a non cooperative game
between incarnations of the decision maker at different instants of time. Nash
equilibrium of these strategies are then considered to define the new concept
of solution of the original problem. Strotz in \cite{Strotz} was the first who
proposed a game theoretic formulation to handle the dynamic time inconsistent
optimal decision problem on the deterministic Ramsey problem, see \cite{Rams}.
Then by capturing the idea of non-commitment, by letting the commitment period
being infinitesimally small, he provided a primitive notion of Nash
equilibrium strategy. Further work along this line in continuous and discrete
time had been done by Pollak \cite{Pollak}, Phelps and Pollak \cite{Phelps},
Goldman \cite{Goldman}, Barro \cite{Barro} and Krusell \& Smith \cite{Krusell}%
. Keeping the same game theoretic approach, Ekland \& Lazrak \cite{EkeLaz} and
Mar\'{\i}n-Solano \& Navas \cite{Solano} treated the optimal consumption
problem where the utility involves a non-exponential discount function in the
deterministic framework. They characterized the equilibrium strategies by a
value function which must satisfy a certain "extended HJB equation", which is
a non linear differential equation displaying a non local term, a term which
depends on the global behaviour of the solution. In this situation, every
decision at time $t$ is taken by a $t-$agent which represents the incarnation
of the controller at time $t$ and is referred in \cite{Solano} as a
"sophisticated $t-$agent".

Bj\"{o}rk \& Murgoci in \cite{BjorkMurgoci} extends the idea to the stochastic
setting where the controlled dynamic is driven by a quite general class of
Markov process and a fairly general objective function. Yong in
\cite{Yong2011}, by a discretization of time, studied a class of time
inconsistent deterministic linear quadratic models and derive equilibrium
controls via some class of Riccati-Voltera equations. Yong in \cite{Yong2012},
also by a discretization of time, investigated a general discounting time
inconsistent stochastic optimal control problem and characterizes a feedback
time-consistent Nash equilibrium control via the so-called "equilibrium HJB
equation". In a series of papers, Basak \& Chabakauri \cite{BasakChaba}, Hu et
al. \cite{Huetal2011}, Czichowsky \cite{Czich} and Bj\"{o}rk et al.
\cite{BjorkMurgociZhou} look at the mean variance problem which is also time inconsistent.

Concerning equilibrium strategies for an optimal consumption-investment
problem with a general discount function, Ekeland \& Pirvu\ \cite{EkePir} are
the first to investigate Nash equilibrium strategies where the price process
of the risky asset is driven by geometric Brownian motion. They characterize
the equilibrium strategies through the solutions of a flow of BSDEs, and they
show, for an special form of the discount function, that the BSDEs reduce to a
system of two ODEs which has a solution. Ekeland et al. in \cite{Ekelandetal}
added life insurance to the investor's portfolio and they characterize the
equilibrium strategy by an integral equation. In \cite{Yong2012}, Yong
discussed the case of time-inconsistent consumption-investment problem under a
power utility function. Following Yong's approach, Zhao et al., in
\cite{Zhaoaetal}, studied the consumption-investment problem with a general
discount function and a logarithmic utility function. Recently, Zou et al., in
\cite{Zoua}, investigated equilibrium consumption-investment decisions for
Merton's portfolio problem with stochastic hyperbolic discounting.\bigskip

\textit{Novelty and Contribution}\bigskip

The purpose of this paper is to investigate equilibrium solutions for a
time-inconsistent consumption-investment problem with a non-exponential
discount function and a general utility function. Different from \cite{Solano}
and \cite{EkePir} where the authors derived explicit solutions for special
forms of the discount factor, in our model, the non-exponential discount
function is in a fairly general form. Moreover, we consider equilibrium
strategies in the open-loop sense, as defined in \cite{Huetal2011} and
\cite{Huetal2015}, which is different from most of the existing literature on
this topic. Note also that the time-inconsistency, in our paper, arises from a
non exponential discounting in the objective function, while the works
\cite{Huetal2011} and \cite{Huetal2015} are concerned with a quite different
kind of time-inconsistency which is caused by the presence of non linear term
of expectations in the terminal cost. On other hand, the objective functional,
in our paper, is not reduced to the quadratic form as in \cite{Huetal2011} and
\cite{Huetal2015}.

By imposing standard assumptions of the classical stochastic maximum
principle, we focus on a variational technique approach leading to a version
of necessary and sufficient condition for equilibrium, which involves a flow
of forward-backward stochastic differential equations (FBSDEs) along with a
certain equilibrium condition. We also present a verification theorem that
covers some possible examples of utility functions. Then, by decoupling the
flow of the FBSDEs, we derive a closed-loop representation of the equilibrium
strategies via some parabolic non-linear partial differential equation (PDE).
We show that within a special form of the utility function (logarithmic, power
and exponential) the PDE reduces to a system of ODEs which has an explicit solution.

We accentuate that, different from most of the existing literature on this
topic, where some feedback equilibrium strategies are derived via several very
complicated highly non-linear integro-differential equations, an explicit
representation of the equilibrium strategies are obtained in our work via
simple ODEs. In addition, this method can provide the necessary and sufficient
conditions to characterize the equilibrium strategies, while the extended HJB
techniques can create, in general, only the sufficient condition in the form
of a verification theorem that characterizes the equilibrium
strategies.\bigskip

\textit{Structure of the paper}\bigskip

The rest of the paper is organized as follows. In Section 2, we formulate the
problem and give the necessary notations and preliminaries. In Section 3 we
present the main results of the paper, Theorem 3.2 and Theorem 3.5 that
characterizes the equilibrium decisions by some necessary and sufficient
conditions. In Section 4, we derive an explicit representation of the
equilibrium consumption-investment strategy. Section 5 is devoted to some
comparisons with existing results in the literature. The paper ends with an
Appendix containing some proofs.

\section{Problem formulation}

Throughout this paper, $(\Omega, \mathcal{F}, {\mathbb{F}},\mathbb{P})$ will
be a filtered probability space such that ${\mathbb{F}}:=\left(
\mathcal{F}_{t}\right)  _{t\in\left[  0,T\right]  }$ is a filtration that
satisfies the usual conditions, in particular, $\mathcal{F}_{0}$ contains all
$\mathbb{P}$-null sets and $\mathcal{F}_{T}=\mathcal{F}$ for an arbitrarily
fixed finite time horizon $T>0.$ Recall that $\mathcal{F}_{t}$ stands for the
information available up to time $t$ and any decision made at time $t$ is
based on this information.

We assume that all processes and random variables are well defined and adapted
to this filtered probability space. In particular, a $d-$dimensional standard
Brownian motion%
\[
W\left(  \cdot\right)  =\left(  W_{1}\left(  \cdot\right)  ,\dots,W_{d}\left(
\cdot\right)  \right)  .
\]
is defined on $(\Omega,\mathcal{F},{\mathbb{F}},\mathbb{P}).$

\subsection{Notations}

Throughout this paper, we use the following notations: $M^{\top}$: the
transpose of the vector (or matrix) $M$, $\left\langle \chi,\zeta\right\rangle
:$ the inner product of $\chi$ and $\zeta$, that is, $\left\langle \chi
,\zeta\right\rangle :=tr(\chi^{T}\zeta).$ For a function $f,$ we denote by
$f_{x}$ (resp. $f_{xx})$ the first (resp. the second) derivative of $f$ with
respect to the variable $x.$

For any Euclidean space $E$ with Frobenius norm $\left\vert \cdot\right\vert $
we let, for any $t\in\left[  0,T\right]  ,$

\begin{enumerate}
\item[$\bullet$] $\mathbb{L}^{p}\left(  \Omega,\mathcal{F}_{t},\mathbb{P}%
;E\right)  :$ for any $p\geq1,$ the set of $E-$valued $\mathcal{F}_{t}%
-$measurable random variables $X,$ such that $\mathbb{E}\left[  \left\vert
X\right\vert ^{p}\right]  <\infty.$

\item[$\bullet$] $\mathcal{L}_{\mathcal{F}}^{1}\left(  t,T;E\right)  :$ the
space of $E-$valued, $\left(  \mathcal{F}_{s}\right)  _{s\in\left[
t,T\right]  }-$adapted processes $c\left(  \cdot\right)  $, with%
\[
\left\Vert c\left(  \cdot\right)  \right\Vert _{\mathcal{L}_{\mathcal{F}}%
^{1}\left(  t,T;E\right)  }=\mathbb{E}\left[  \int_{0}^{T}\left\vert c\left(
t\right)  \right\vert dt\right]  <\infty.
\]

\item[$\bullet$] $\mathcal{L}_{\mathcal{F}}^{2}\left(  t,T;E\right)  :$ the
space of $E-$valued, $\left(  \mathcal{F}_{s}\right)  _{s\in\left[
t,T\right]  }-$adapted continuous processes $Y\left(  \cdot\right)  $, with%
\[
\ \left\Vert Y\left(  \cdot\right)  \right\Vert _{\mathcal{L}_{\mathcal{F}%
}^{2}\left(  t,T;E\right)  }=\sqrt{\mathbb{E}\left[  \sup\limits_{s\in\left[
t,T\right]  }\left\vert Y\left(  s\right)  \right\vert ^{2}\right]  }<\infty.
\]

\item[$\bullet$] $\mathcal{M}_{\mathcal{F}}^{2}\left(  t,T;E\right)  :$ the
space of $E-$valued, $\left(  \mathcal{F}_{s}\right)  _{s\in\left[
t,T\right]  }-$adapted processes $Z\left(  \cdot\right)  $, with%
\[
\ \left\Vert Z\left(  \cdot\right)  \right\Vert _{\mathcal{M}_{\mathcal{F}%
}^{2}\left(  t,T;E\right)  }=\sqrt{\mathbb{E}\left[  {\displaystyle\int
_{t}^{T}}\left\vert Z\left(  s\right)  \right\vert ^{2}ds\right]  }<\infty.
\]

\end{enumerate}

\subsection{Financial market}

Consider an individual facing the inter-temporal consumption and portfolio
problem where the market environment consists of one riskless and $d$ risky
securities. The risky securities are stocks and their prices are modeled as
It\^{o} processes. Namely, for $i=1,2,..,d,$ the price $S_{i}\left(  s\right)
,$ for $s\in\left[  0,T\right]  ,$ of the i-th risky asset, satisfies%

\begin{equation}
dS_{i}\left(  s\right)  =S_{i}\left(  s\right)  \left(  \mu_{i}\left(
s\right)  ds +\sum_{j=1}^{d}\sigma_{ij}\left(  s\right)  dW_{j}\left(
s\right)  \right)  ,\text{ } \tag{2.1}%
\end{equation}
with $S_{i}\left(  0\right)  >0,$ for $i=1,2,...,d,$ and the coefficients
$\mu_{i}\left(  \cdot\right)  $ and $\sigma_{i}\left(  \cdot\right)  =\left(
\sigma_{i1}\left(  \cdot\right)  ,\dots,\sigma_{id}\left(  \cdot\right)
\right)  ,$ for $i=1,..,d,$ are ${\mathbb{F}}-$progressively measurable
processes with values in $\mathbb{R}$ and $\mathbb{R}^{d}$, respectively. For
brevity, we use $\mu\left(  \cdot\right)  =\left(  \mu_{1}\left(
\cdot\right)  ,\mu_{2}\left(  \cdot\right)  ,\dots, \mu_{d}\left(
\cdot\right)  \right)  $ to denote the drift rate vector and $\sigma\left(
\cdot\right)  = \left(  \sigma_{ij}\left(  \cdot\right)  \right)  _{1\leq
i,j\leq d}$ to denote the random volatility matrix.

The riskless asset, or the savings account, has the price process
$S_{0}\left(  s\right)  $, for $s\in\left[  0,T\right]  ,$ governed by%

\begin{equation}
dS_{0}\left(  s\right)  =r_{0}\left(  s\right)  S_{0}\left(  s\right)
ds,\text{ }S_{0}\left(  0\right)  =1, \tag{2.2}%
\end{equation}
where $r_{0}\left(  \cdot\right)  $ is a deterministic function with values in
$\left[  0,\infty\right)  $ that represents the interest rate. We assume that
$\mathbb{E}\left[  \mu_{i}\left(  t\right)  \right]  >r_{0}\left(  t\right)
\geq0,$ $dt-a.e.,$ for $i=1,2,..,d$. This is a very natural assumption, since
otherwise, nobody is willing to invest in the risky stocks.

\subsection{Investment-consumption policies and wealth process}

Starting from an initial capital $x_{0}>0$ at time $0$, during the time
horizon $\left[  0,T\right]  $, the decision maker is allowed to dynamically
invest in the stocks as well as in the bond and consuming. A
consumption-investment strategy is described by a $\left(  d+1\right)
$-dimensional stochastic process $u\left(  \cdot\right)  =\left(  c\left(
\cdot\right)  ,u_{1}\left(  \cdot\right)  ,\dots, u_{d}\left(  \cdot\right)
\right)  ^{\top},$ where $c\left(  s\right)  $ represents the consumption rate
at time $s\in\left[  0,T\right]  $ and $u_{i}\left(  s\right)  ,$ for
$i=1,2,..,d,$ represents the amount invested in the $i$-$th$ risky stock at
time $s\in\left[  0,T\right]  .$ The process $u_{I}\left(  \cdot\right)
=\left(  u_{1}\left(  \cdot\right)  ,\dots, u_{d}\left(  \cdot\right)
\right)  ^{\top}$ is called an investment strategy. The amount invested in the
bond at time $s$ is%

\[
X^{x_{0}, u}\left(  s\right)  -\sum\limits_{i=1}^{d}u_{i}\left(  s\right)  ,
\]
where $X^{x_{0},u}\left(  \cdot\right)  $ is the wealth process associated
with the strategy $u\left(  \cdot\right)  $ and the initial capital $x_{0}$.
The evolution of $X^{x_{0},u}\left(  \cdot\right)  $ can be described as%

\[
\left\{
\begin{array}
[c]{l}%
dX^{x_{0},u}\left(  s\right)  =\left(  X^{x_{0},u}\left(  s\right)
-\sum\limits_{i=1}^{d}u_{i}\left(  s\right)  \right)  \dfrac{dS_{0}\left(
s\right)  }{S_{0}\left(  s\right)  } +\sum\limits_{i=1}^{d}u_{i}\left(
s\right)  \dfrac{dS_{i}\left(  s\right)  }{S_{i}\left(  s\right)  }-c\left(
s\right)  ds, \text{ for }s\in\left[  0,T\right]  ,\\
X^{x_{0},u}\left(  0\right)  =x_{0}.
\end{array}
\right.
\]

Accordingly, the wealth process solves the SDE%

\begin{equation}
\left\{
\begin{array}
[c]{l}%
dX^{x_{0},u}\left(  s\right)  =\left\{  r_{0}\left(  s\right)  X^{x_{0}%
,u}\left(  s\right)  +u_{I}\left(  s\right)  ^{\top} r\left(  s\right)
-c\left(  s\right)  \right\}  ds\\
\text{ \ \ \ \ \ \ \ \ \ \ \ \ \ \ \ \ \ \ \ \ } +u_{I}\left(  s\right)
^{\top}\sigma\left(  s\right)  dW\left(  s\right)  ,\text{ for }s\in\left[
0,T\right]  ,\\
X^{x_{0},u}\left(  0\right)  =x_{0}.
\end{array}
\right.  \tag{2.3}%
\end{equation}
where $r\left(  \cdot\right)  =\left(  \mu_{1}\left(  \cdot\right)  -r_{0}
\left(  \cdot\right)  ,\dots, \mu_{d}\left(  \cdot\right)  -r_{0}\left(
\cdot\right)  \right)  ^{\top}$.

As time evolves, it is natural to consider the controlled stochastic
differential equation parametrized by $\left(  t,\xi\right)  \in\left[
0,T\right]  \times\mathbb{L}^{2}\left(  \Omega,\mathcal{F}_{t},\mathbb{P}%
;\mathbb{R}\right)  $ and satisfied by $X\left(  \cdot\right)  =X^{t,\xi
}\left(  \cdot;u\left(  \cdot\right)  \right)  ,$%
\begin{equation}
\left\{
\begin{array}
[c]{l}%
dX\left(  s\right)  =\left\{  r_{0}\left(  s\right)  X\left(  s\right)
+u_{I}\left(  s\right)  ^{\top}r\left(  s\right)  -c\left(  s\right)
\right\}  ds+u_{I}\left(  s\right)  ^{\top}\sigma\left(  s\right)  dW\left(
s\right)  ,\text{ for }s\in\left[  t,T\right]  ,\\
X\left(  t\right)  =\xi.
\end{array}
\right.  \tag{2.4}%
\end{equation}

\begin{definition}
[Admissible Strategy]A strategy $u\left(  \cdot\right)  =\left(  c\left(
\cdot\right)  ,u_{I}\left(  \cdot\right)  ^{\top}\right)  ^{\top}$ is said to
be admissible over $\left[  t,T\right]  $ if $u\left(  \cdot\right)
\in\mathcal{L}_{\mathcal{F}}^{1}\left(  t,T;\mathbb{R}\right)  \times
\mathcal{M}_{\mathcal{F}}^{2}\left(  t,T;\mathbb{R}^{d}\right)  $ and for any
$\left(  t,\xi\right)  \in\left[  0,T\right]  \times\mathbb{L}^{2}\left(
\Omega,\mathcal{F}_{t},\mathbb{P};\mathbb{R}\right)  ,$ the equation $\left(
2.4\right)  $ has a unique solution $X\left(  \cdot\right)  =X^{t,\xi}\left(
\cdot;u\left(  \cdot\right)  \right)  .$
\end{definition}

We impose the following assumption about the coefficients.

\begin{enumerate}
\item[\textbf{(H1)}] Processes $r_{0}\left(  \cdot\right)  $, $r\left(
\cdot\right)  $ and $\sigma\left(  \cdot\right)  $ are uniformly bounded.
Moreover we assume the following uniform ellipticity condition:%

\[
\sigma\left(  s\right)  \sigma\left(  s\right)  ^{\top}\geq\epsilon
I_{d},\text{ }ds-a.e,\text{ }d\mathbb{P-}a.s.
\]
for some $\epsilon>0$, where $I_{d}$ denotes the identity matrix on
$\mathbb{R}^{d\times d}.$
\end{enumerate}

Under (\textbf{H1)}, for any $\left(  t,\xi,u\left(  \cdot\right)  \right)
\in\left[  0,T\right]  \times\mathbb{L}^{2}\left(  \Omega,\mathcal{F}%
_{t},\mathbb{P};\mathbb{R}\right)  \times\mathcal{L}_{\mathcal{F}}^{1}\left(
t,T;\mathbb{R}\right)  \times\mathcal{M}_{\mathcal{F}}^{2}\left(
t,T;\mathbb{R}^{d}\right)  ,$ the state equation $\left(  2.4\right)  $ has a
unique solution $X\left(  \cdot\right)  \in\mathcal{L}_{\mathcal{F}}%
^{2}\left(  t,T;\mathbb{R}\right)  $. Moreover, we have the estimate%

\begin{equation}
\mathbb{E}\left[  \sup_{t\leq s\leq T}\left\vert X\left(  s\right)
\right\vert ^{2}\right]  \leq K\left(  1+\mathbb{E}\left[  \left\vert
\xi\right\vert ^{2}\right]  \right)  , \tag{2.5}%
\end{equation}
for some positive constant $K$. In particular for $t=0,$ $x_{0}>0$ and
$u\left(  \cdot\right)  =\left(  c\left(  \cdot\right)  ,u_{I}\left(
\cdot\right)  ^{\top}\right)  ^{\top}\in\mathcal{L}_{\mathcal{F}}^{1}\left(
0,T;\mathbb{R}\right)  \times\mathcal{M}_{\mathcal{F}}^{2}\left(
0,T;\mathbb{R}^{d}\right)  ,$ the state equation $\left(  2.3\right)  $ has a
unique solution $X^{x_{0},u}\left(  \cdot\right)  \in\mathcal{L}_{\mathcal{F}%
}^{2}\left(  0,T;\mathbb{R}\right)  $ and the following estimate holds:%

\begin{equation}
\mathbb{E}\left[  \sup_{0\leq s\leq T}\left\vert X^{x_{0},u}\left(  s\right)
\right\vert ^{2}\right]  \leq K\left(  1+\left\vert x_{0}\right\vert
^{2}\right)  . \tag{2.6}%
\end{equation}

\subsection{General discounted utility function}

Most of financial-economics works have considered that the rate of time
preference is constant (exponential discounting). However there is growing
evidence to suggest that this may not be the case. In this section, we discuss
the general discounting preferences. We also introduce the basic modeling
framework of Merton's consumption and portfolio problem. We refer the reader
to \cite{Flemin}, \cite{Karat}, \cite{Mert1}, \cite{Mert2} and \cite{Plisk}
for more detail about the classical Merton model.

\subsubsection{Discount function}

As soon as discounting is non-exponential, most papers work with special form
of the non-exponential discount factor. Different to these works we consider a
general form of the discount factor.

\begin{definition}
A discount function $\lambda:\left[  0,T\right]  \rightarrow%
\mathbb{R}
$ is a deterministic function satisfying $\lambda\left(  0\right)  =1,$
$\lambda\left(  s\right)  >0$ $ds-a.e.$ and $\int_{0}^{T}\lambda\left(
s\right)  ds<\infty.$
\end{definition}

We also impose the following Lipschitz condition with constant $C,$ on
$\lambda\left(  .\right)  $

\begin{enumerate}
\item[\textbf{(H2)}] There exists a constant $C>0$ such that $\left\vert
\lambda\left(  s\right)  -\lambda\left(  t\right)  \right\vert \leq
C\left\vert s-t\right\vert ,$ for any $t,s\in\left[  0,T\right]  .$
\end{enumerate}

\begin{remark}
Note that Assumption (\textbf{H2)} is satisfied by many discount functions,
such as exponential discount functions, see \cite{Mert1} and \cite{Mert2},
mixture of exponential functions, see \cite{EkePir}, and hyperbolic discount
functions, see \cite{Zhaoaetal}.
\end{remark}

\subsubsection{Utility functions and objective}

In order to evaluate the performance of a consumption-investment strategy, the
decision maker derives utility from inter-temporal consumption and final
wealth. Let $\varphi\left(  \cdot\right)  $ be the utility of inter-temporal
consumption and $h\left(  \cdot\right)  $ the utility of the terminal wealth
at some non-random horizon $T$ (which is a primitive of the model). Then, for
any $\left(  t,\xi\right)  \in\left[  0,T\right]  \times\mathbb{L}^{2}\left(
\Omega,\mathcal{F}_{t},\mathbb{P};\mathbb{R}\right)  $ the
investment-consumption optimization problem is reduced to maximize the utility
function $J\left(  t,\xi,.\right)  $ given by%
\begin{equation}
J\left(  t,\xi,u\left(  \cdot\right)  \right)  =\mathbb{E}^{t}\left[  \int
_{t}^{T}\lambda\left(  s-t\right)  \varphi\left(  c\left(  s\right)  \right)
ds+\lambda\left(  T-t\right)  h\left(  X\left(  T\right)  \right)  \right]  ,
\tag{2.7}%
\end{equation}
over $u\left(  \cdot\right)  \in\mathcal{L}_{\mathcal{F}}^{1}\left(
t,T;\mathbb{R}\right)  \times\mathcal{M}_{\mathcal{F}}^{2}\left(
t,T;\mathbb{R}^{d}\right)  ,$ subject to $\left(  2.4\right)  ,$ where
$\mathbb{E}^{t}\left[  \cdot\right]  =\mathbb{E}\left[  \cdot\left\vert
\mathcal{F}_{t}\right.  \right]  $. We restrict ourselves to utility functions
which satisfy the following condition:

\begin{enumerate}
\item[\textbf{(H2)}] The maps $\varphi\left(  \cdot\right)  ,h\left(
\cdot\right)  :[0,\infty)\rightarrow\mathbb{R}$ are strictly increasing,
strictly concave and satisfy the Inada conditions.
\end{enumerate}

If we write $W^{\star}\left(  s\right)  =\left(  0,W\left(  s\right)  ^{\top
}\right)  ^{\top}$ and we denote $B\left(  s\right)  =\left(  -1,r\left(
s\right)  ^{\top}\right)  ^{\top},$ $\Gamma=\left(  1,0_{\mathbb{R}^{d}}%
^{\top}\right)  ^{\top}$ and%

\[
D\left(  s\right)  =\left(
\begin{array}
[c]{cc}%
0 & 0_{\mathbb{R}^{d}}^{\top}\\
0_{\mathbb{R}^{d}} & \sigma\left(  s\right)
\end{array}
\right)  ,
\]
then the optimal control problem associated with $\left(  2.4\right)  $ and
$\left(  2.7\right)  $ is equivalent to maximize%

\begin{equation}
J\left(  t,\xi,u\left(  \cdot\right)  \right)  =\mathbb{E}^{t}\left[  \int
_{t}^{T}\lambda\left(  s-t\right)  \varphi\left(  \Gamma^{\top}u\left(
\cdot\right)  \right)  ds +\lambda\left(  T-t\right)  h\left(  X\left(
T\right)  \right)  \right]  , \tag{2.8}%
\end{equation}
subject to%

\begin{equation}
\left\{
\begin{array}
[c]{l}%
dX\left(  s\right)  =\left\{  r_{0}\left(  s\right)  X\left(  s\right)
+u\left(  s\right)  ^{\top}B\left(  s\right)  \right\}  ds+u\left(  s\right)
^{\top}D\left(  s\right)  dW^{\star}\left(  s\right)  , \text{ for }%
s\in\left[  t,T\right]  ,\\
X\left(  t\right)  =\xi.
\end{array}
\right.  \tag{2.9}%
\end{equation}

\subsubsection{Time inconsistency}

Let us first note that the optimal policies, although they exist, will not be
time-consistent in general. First of all, as an illustration, let us consider
the model in (2.8)--(2.9) with logarithmic utility functions. We suppose that
the financial market consists of one riskless asset and $d$ risky assets.
Arguing as in \cite{EkePir}$,$ we can prove that, if the agent is naive and
starts with a given positive wealth $x$, at some instant $t,$ then by the
standard dynamic programming approach, the value function associated with this
stochastic control problem solves the following Hamilton--Jacobi--Bellman equation%

\begin{equation}
\left\{
\begin{array}
[c]{l}%
V_{s}^{t}\left(  s,x\right)  +\underset{\left(  c,u_{I}\right)  \in
\mathbb{R}^{d+1}}{\sup} \left\{  \left(  r_{0}\left(  s\right)  X\left(
s\right)  +u_{I}{}^{\top}r\left(  s\right)  -c\right)  V_{x}^{t}\left(
s,x\right)  \right.  +\dfrac{1}{2}u_{I}^{\top}\sigma\left(  s\right)
\sigma\left(  s\right)  ^{\top}u_{I}V_{xx}^{t}\left(  s,x\right) \\
\text{ \ \ \ \ \ \ \ \ \ \ \ \ \ \ \ \ \ \ \ \ \ \ \ \ \ \ \ }\left.
+\dfrac{\lambda^{\prime}\left(  s-t\right)  }{\lambda\left(  s-t\right)
}V^{t}\left(  s,x\right)  +\varphi\left(  c\right)  \right\}  =0,\bigskip
\text{ for }s\in\left[  t,T\right]  ,\\
V^{t}\left(  T,x\right)  =h\left(  x\right)  .
\end{array}
\right.  \tag{2.10}%
\end{equation}

The HJB equation contains the term $\dfrac{\lambda^{\prime}\left(  s-t\right)
}{\lambda\left(  s-t\right)  }$, which depends not only on the current time
$s$ but also on initial time $t$, and so, the optimal policy will depend on
$t$ as well. Indeed, the first order necessary conditions yield the
$t-$optimal policy%

\begin{align*}
\overline{u}_{I}^{t}\left(  s\right)   &  =r\left(  s\right)  \left(
\sigma\left(  s\right)  \sigma\left(  s\right)  ^{\top}\right)  ^{-1}%
\frac{V_{x}^{t}\left(  s,x\right)  }{V_{xx}^{t}\left(  s,x\right)  },\\
\overline{c}^{t}\left(  s\right)   &  =\varphi^{-1}\left(  V_{x}^{t}\left(
s,x\right)  \right)  .
\end{align*}

Let us consider the following example: $\varphi\left(  x\right)  =h\left(
x\right)  =\log x.$ The naive agent for the initial pair $\left(
0,x_{0}\right)  $ solves the problem, assuming that the discount rate of time
preference will be $\lambda\left(  s\right)  $, for $s\in\left[  0,T\right]
,$ and the optimal consumption strategy will be%

\[
\overline{c}^{0,x_{0}}\left(  s\right)  =\left[  1+\int_{s}^{T} \exp\left\{
\lambda\left(  r-s\right)  +\log\left(  \frac{\lambda\left(  r\right)
}{\lambda\left(  s\right)  }\right)  \right\}  dr\right]  ^{-1},\text{ for }
s\in\left[  0,T \right]  .
\]

This solution corresponds to the so-called pre-commitment solution, in the
sense that it is optimal as long as the agent can precommit (by signing a
contract, for example) his or her future behavior at time $t=0$. If there is
no commitment, the 0-agent will take the action $\overline{c}^{0,x_{0}}\left(
s\right)  $ but, in the near future, the $\epsilon$-agent will change his
decision rule (time-inconsistency) to the solution of the HJB equation
$\left(  2.10\right)  $ with $t=\epsilon.$ In this case the optimal control
trajectory for $s>\epsilon$ will be changed to $\overline{c}^{\epsilon
,x_{\epsilon}}\left(  s\right)  $ given by%
\[
\overline{c}^{\epsilon,x_{\epsilon}}\left(  s\right)  =\overline{c}%
^{\epsilon,\bar{X}\left(  \epsilon\right)  }\left(  s\right)  =\left[
1+\int_{s}^{T}\exp\left\{  \lambda\left(  r-s\right)  +\log\left(
\frac{\lambda\left(  r-\epsilon\right)  }{\lambda\left(  s-\epsilon\right)
}\right)  \right\}  dr\right]  ^{-1},\text{ for }s\in\left[  \epsilon
,T\right]  .
\]

If $\lambda\left(  t\right)  =e^{-\delta t}$ where $\delta>0$ is the
constant\textbf{ }discount rate, then
\[
\overline{c}_{|\left[  \epsilon,T\right]  }^{0,x_{0}}\left(  s\right)
=\overline{c}^{\epsilon,x_{\epsilon}}\left(  s\right)  ,\text{ for }%
s\in\left[  \epsilon,T\right]  ,
\]
hence the optimal consumption plan is time consistent. As soon as discount
function is non-exponential%

\[
\overline{c}_{|\left[  \epsilon,T\right]  }^{0,x_{0}}\left(  s\right)
\neq\overline{c}^{\epsilon,x_{\epsilon}}\left(  s\right)  ,\text{ for }
s\in\left[  \epsilon,T\right]  .
\]

Then the optimal consumption plan is not time consistent. In general, the
solution for the naive agent will be constructed by solving the family of HJB
equations $\left(  2.10\right)  $ for $t\in\left[  0,T\right]  $, and patching
together the \textquotedblleft optimal\textquotedblright\ solutions
$\overline{c}^{t,x_{t}}\left(  t\right)  .$ If the agent is sophisticated,
things become more complicated. The standard HJB equation cannot be used to
construct the solution, and a new method is required in what follows.

\section{Equilibrium strategies}

It is well known that the problem described above by $\left(  2.8\right)
-\left(  2.9\right)  $ turns out to be time inconsistent in the sense that it
does not satisfy the Bellman optimality principle, since a restriction of an
optimal control for a specific initial pair on a later time interval might not
be optimal for that corresponding initial pair. For a more detailed discussion
see Ekeland \& Pirvu\ \cite{EkePir} and Yong \cite{Yong2012}. Since the lack
of time consistency, we consider open-loop Nash equilibrium controls instead
of optimal controls. As in \cite{Huetal2011}, we first consider an equilibrium
by local spike variation, given, for $t\in\left[  0,T\right]  ,$ an admissible
consumption-investment strategy $\hat{u}\left(  \cdot\right)  \in
\mathcal{L}_{\mathcal{F}}^{1}\left(  t,T;\mathbb{R}\right)  \times
\mathcal{M}_{\mathcal{F}}^{2}\left(  t,T;\mathbb{R}^{d}\right)  .$ For any
$\mathbb{R}^{d+1}-$valued, $\mathcal{F}_{t}-$measurable and bounded random
variable $v$ and for any $\varepsilon>0,$ define%

\begin{equation}
u^{\varepsilon}\left(  s\right)  :=\left\{
\begin{array}
[c]{ll}%
\hat{u}\left(  s\right)  +v, & \text{ for }s\in\left[  t,t+\varepsilon\right)
,\\
\hat{u}\left(  s\right)  , & \text{ for }s\in\left[  t+\varepsilon,T\right]  .
\end{array}
\right.  \tag{3.1}%
\end{equation}

We have the following definition.

\begin{definition}
[Open-loop Nash equilibrium]An admissible strategy $\hat{u}\left(
\cdot\right)  \in\mathcal{L}_{\mathcal{F}}^{1}\left(  t,T;\mathbb{R}\right)
\times\mathcal{M}_{\mathcal{F}}^{2}\left(  t,T;\mathbb{R}^{d}\right)  $ is an
open-loop Nash equilibrium strategy if%

\begin{equation}
\lim_{\varepsilon\downarrow0}\inf\frac{1}{\varepsilon}\left\{  J\left(
t,\hat{X}\left(  t\right)  ,u^{\varepsilon}\left(  \cdot\right)  \right)
-J\left(  t,\hat{X}\left(  t\right)  ,\hat{u}\left(  \cdot\right)  \right)
\right\}  \leq0, \tag{3.2}%
\end{equation}
for any $t\in\left[  0,T\right]  ,$ where $\hat{X}$ is the equilibrium wealth
process that solves the SDE%

\begin{equation}
\left\{
\begin{array}
[c]{l}%
d\hat{X}\left(  s\right)  =\left\{  r_{0}\left(  s\right)  \hat{X}\left(
s\right)  +\hat{u}\left(  s\right)  ^{\top}B\left(  s\right)  \right\}
ds+\hat{u}\left(  s\right)  ^{\top}D\left(  s\right)  dW^{\star}\left(
s\right)  ,\text{ for }s\in\left[  t,T\right]  ,\\
\hat{X}\left(  t\right)  =\xi.
\end{array}
\right.  \tag{3.3}%
\end{equation}

\end{definition}

\subsection{A necessary and sufficient condition for equilibrium controls}

In this paper we follow an alternative approach, which is essentially a
necessary and sufficient conditions for equilibrium. In the same spirit of
proving the stochastic Pontryagin's maximum principle for equilibrium in
\cite{Huetal2011} for the case of linear quadratic models, we derive this
condition by a second-order expansion in the spike variation. First, during
this section we impose the following hypothesis on the utility functions

\begin{enumerate}
\item[(\textbf{H3})] The maps $\varphi\left(  \cdot\right)  ,h\left(
\cdot\right)  $ are twice continuously differentiable functions. We suppose
also that, there exists a positive constant $C$ such that%
\[
\left\vert \varphi_{xx}\left(  x\right)  -\varphi_{xx}\left(  \hat{x}\right)
\right\vert +\left\vert h_{xx}\left(  x\right)  -h_{xx}\left(  \hat{x}\right)
\right\vert \leq C\left\vert x-\hat{x}\right\vert ,\forall x,\hat{x}\in
\lbrack0,\infty).
\]

\end{enumerate}

Now, we introduce the adjoint equations involved in the characterization of
open-loop Nash equilibrium controls.

\subsubsection{Adjoint processes}

Let $\hat{u}\left(  \cdot\right)  =\left(  \hat{c}\left(  \cdot\right)
,\hat{u}_{I}\left(  \cdot\right)  ^{\top}\right)  ^{\top}\in\mathcal{L}%
_{\mathcal{F}}^{1}\left(  0,T;\mathbb{R}\right)  \times\mathcal{M}%
_{\mathcal{F}}^{2}\left(  0,T;\mathbb{R}^{d}\right)  $ an admissible strategy
and denote by $\hat{X}\left(  \cdot\right)  \in\mathcal{L}_{\mathcal{F}}%
^{2}\left(  0,T;\mathbb{R}\right)  $ the corresponding wealth process. For
each $t\in\left[  0,T\right]  $, we introduce the first order adjoint equation
defined on the time interval $\left[  t,T\right]  $, and satisfied by the pair
of processes $\left(  p\left(  \cdot;t\right)  ,q\left(  \cdot;t\right)
\right)  $ as follows%

\begin{equation}
\left\{
\begin{array}
[c]{l}%
dp\left(  s;t\right)  =-r_{0}\left(  s\right)  p\left(  s;t\right)
ds+q\left(  s;t\right)  ^{\top}dW\left(  s\right)  ,\text{ for }s\in\left[
t,T\right]  ,\\
p\left(  T;t\right)  =\lambda\left(  T-t\right)  h_{x}\left(  \hat{X}\left(
T\right)  \right)  ,
\end{array}
\right.  \tag{3.4}%
\end{equation}
where $q\left(  \cdot;t\right)  =\left(  q_{1}\left(  \cdot;t\right)  ,\dots,
q_{d}\left(  \cdot;t\right)  \right)  ^{\top}.$ Under the assumption
(\textbf{H1),} the BSDE $\left(  3.4\right)  $ is uniquely solvable in
$\mathcal{L}_{\mathcal{F}}^{2}\left(  t,T;\mathbb{R}\right)  \times
\mathcal{M}_{\mathcal{F}}^{2}\left(  0,T;\mathbb{R}^{d}\right)  .$ Moreover
there exists a constant $K>0$ such that, for any $t\in\left[  0,T\right]  ,$
we have the following estimate%

\begin{equation}
\left\Vert p\left(  \cdot;t\right)  \right\Vert _{\mathcal{L}_{\mathcal{F}%
}^{2} \left(  t,T;\mathbb{R}\right)  }^{2}+\left\Vert q\left(  \cdot;t\right)
\right\Vert _{\mathcal{M}^{2}\left(  t,T;\mathbb{R}^{d}\right)  }^{2}\leq
K\left(  1+\xi^{2}\right)  . \tag{3.5}%
\end{equation}

The second order adjoint equation is defined on the time interval $\left[
t,T\right]  $ and satisfied by the pair of processes $\left(  P\left(
\cdot;t\right)  , Q\left(  \cdot;t\right)  \right)  \in\mathcal{L}%
_{\mathcal{F}}^{2}\left(  t,T;\mathbb{R}\right)  \times\mathcal{M}%
_{\mathcal{F}}^{2}\left(  t,T;\mathbb{R}^{d}\right)  $ as follows%
\begin{equation}
\left\{
\begin{array}
[c]{l}%
dP\left(  s;t\right)  =-2r_{0}\left(  s\right)  P\left(  s;t\right)
ds+Q\left(  s;t\right)  ^{\top}dW\left(  s\right)  ,\text{ for }s\in\left[
t,T\right]  ,\\
P\left(  T;t\right)  =\lambda\left(  T-t\right)  h_{xx}\left(  \hat{X}\left(
T\right)  \right)  .
\end{array}
\right.  \tag{3.6}%
\end{equation}
where $Q\left(  \cdot;t\right)  =\left(  Q_{1}\left(  \cdot;t\right)  ,\dots,
Q_{d}\left(  \cdot;t\right)  \right)  ^{\top}.$ Under \textbf{(H1)} the above
BSDE has a unique solution $\left(  P\left(  \cdot;t\right)  , Q\left(
\cdot;t\right)  \right)  \in\mathcal{L}_{\mathcal{F}}^{2}\left(
t,T;\mathbb{R}\right)  \times\mathcal{M}_{\mathcal{F}}^{2} \left(
t,T;\mathbb{R}^{d}\right)  $. Moreover we have the following representation
for $P\left(  \cdot;t\right)  :$%

\begin{equation}%
\begin{array}
[c]{l}%
P\left(  s;t\right)  =\mathbb{E}^{s}\left[  \lambda\left(  T-t\right)
h_{xx}\left(  \hat{X}\left(  T\right)  \right)  e^{\int_{s}^{T}2r_{0}\left(
\tau\right)  d\tau}\right]  ,\text{ for }s\in\left[  t,T\right]  .
\end{array}
\tag{3.7}%
\end{equation}

Indeed, if we define the function $\Theta\left(  \cdot,t\right)  ,$ for each
$t\in\left[  0,T\right]  ,$ as the fundamental solution of the linear ODE%

\begin{equation}
\left\{
\begin{array}
[c]{l}%
d\Theta\left(  \tau,t\right)  =r_{0}\left(  \tau\right)  \Theta\left(
\tau,t\right)  d\tau,\text{ for }\tau\in\left[  t,T\right]  ,\\
\Theta\left(  t,t\right)  =1,
\end{array}
\right.  \tag{3.8}%
\end{equation}
and we apply the It\^{o}'s formula to $\tau\rightarrow P\left(  \tau;t\right)
\Theta\left(  \tau,t \right)  ^{2}$ on $\left[  t,T\right]  ,$ by taking
conditional expectations, we obtain $\left(  3.7\right)  $. Note that since
$h_{xx}\left(  \hat{X}\left(  T\right)  \right)  \leq0$, then $P\left(
s;t\right)  \leq0,$ $ds-a.e$.

\subsubsection{A characterization of equilibrium strategies}

The following theorem is the first main result of this work, it provides a
necessary and sufficient condition for equilibrium. As we have said before,
the proof is inspired in \cite{Huetal2011} and \cite{Huetal2015}.

First, we define the process $\tilde{q}\left(  s;t\right)  =\left(  0,q\left(
s;t\right)  ^{\top}\right)  ^{\top},$ and we introduce the following notations:%

\begin{equation}
\mathcal{H}\left(  s;t\right)  \triangleq p\left(  s;t\right)  B\left(
s\right)  +D\left(  s\right)  \tilde{q}\left(  s;t\right)  +\lambda\left(
s-t\right)  \varphi_{c}\left(  \Gamma^{\top}\hat{u}\left(  s\right)  \right)
\Gamma\tag{3.9}%
\end{equation}
and%

\begin{equation}
\mathcal{A}\left(  s;t\right)  \triangleq\left(
\begin{array}
[c]{cc}%
\lambda\left(  s-t\right)  \varphi_{cc}\left(  \Gamma^{\top}\left(  \hat
{u}\left(  s\right)  \right)  \right)  & 0_{\mathbb{R}^{n}}^{\top}\\
0_{\mathbb{R}^{n}} & \sigma\left(  s\right)  \sigma\left(  s\right)  ^{\top
}P\left(  s;t\right)
\end{array}
\right)  . \tag{3.10}%
\end{equation}

\begin{theorem}
Let \textbf{(H1)-(H3)} hold. Given an admissible strategy $\hat{u}\left(
\cdot\right)  \in\mathcal{L}_{\mathcal{F}}^{1}\left(  0,T;\mathbb{R}\right)
\times\mathcal{M}_{\mathcal{F}}^{2}\left(  0,T;\mathbb{R}^{d}\right)  $, let
for any $t\in\left[  0,T\right]  ,$ the processus
\[
\left(  p\left(  \cdot;t\right)  ,q\left(  \cdot;t\right)  \right)
\in\mathcal{L}_{\mathcal{F}}^{2}\left(  t,T;\mathbb{R}\right)  \times
\mathcal{M}_{\mathcal{F}}^{2}\left(  t,T;\mathbb{R}^{d}\right)
\]
be the unique solution to the BSDE $\left(  3.4\right)  $. Then, $\hat
{u}\left(  \cdot\right)  $ is an equilibrium consumption-investment strategy,
if and only if, the following condition holds%
\begin{equation}
\mathcal{H}\left(  t;t\right)  =0,\text{ }d\mathbb{P-}a.s.,\text{ }dt-a.e.
\tag{3.11}%
\end{equation}

\end{theorem}

In order to derive the proof of this theorem, let us, first of all, derive
some technical results. First, denote by $\hat{X}^{\varepsilon}\left(
\cdot\right)  $ the solution of the state equation corresponding to
$u^{\varepsilon}\left(  \cdot\right)  $. Since the coefficients of the
controlled state equation are linear, using the standard perturbation
approach, see e.g. \cite{YongZhou}, we have%
\begin{equation}
\hat{X}^{\varepsilon}\left(  s\right)  -\hat{X}\left(  s\right)
=y^{\varepsilon,v}\left(  s\right)  +z^{\varepsilon,v}\left(  s\right)
,\text{ for }s\in\left[  t,T\right]  , \tag{3.12}%
\end{equation}
where for any $\mathbb{R}^{d+1}-$valued, $\mathcal{F}_{t}-$measurable and
bounded random variable $v$ and for any $\varepsilon\in\left[  0,T-t\right)
,$ $y^{\varepsilon,v}\left(  \cdot\right)  $ and $z^{\varepsilon,v}\left(
\cdot\right)  $ solve respectively the following linear stochastic
differential equations:%

\begin{equation}
\left\{
\begin{array}
[c]{l}%
dy^{\varepsilon,v}\left(  s\right)  =r_{0}\left(  s\right)  y^{\varepsilon
,v}\left(  s\right)  ds +v^{\top}D\left(  s\right)  1_{\left[  t,t+\varepsilon
\right)  }\left(  s\right)  dW^{\star}\left(  s\right)  ,\text{\ for }%
s\in\left[  t,T\right]  ,\\
y^{\varepsilon,v}\left(  t\right)  =0,
\end{array}
\right.  \tag{3.13}%
\end{equation}
and%

\begin{equation}
\left\{
\begin{array}
[c]{l}%
dz^{\varepsilon,v}\left(  s\right)  =\left\{  r_{0}\left(  s\right)
z^{\varepsilon,v}\left(  s\right)  +v^{\top}B\left(  s\right)  1_{\left[
t,t+\varepsilon\right)  } \left(  s\right)  \right\}  ds,\text{\ for }%
s\in\left[  t,T\right]  ,\\
z^{\varepsilon,v}\left(  t\right)  =0.
\end{array}
\right.  \tag{3.14}%
\end{equation}

\begin{proposition}
Let \textbf{(H1) }holds. For any $t\in\left[  0,T\right]  ,$ the following
estimates hold for any $k\geq1:$%
\begin{align}
&  \mathbb{E}^{t}\left[  \sup\limits_{s\in\left[  t,T\right]  }\left\vert
y^{\varepsilon,v}\left(  s\right)  \right\vert ^{2k}\right]  =O\left(
\varepsilon^{k}\right)  ,\tag{3.15}\\
&  \mathbb{E}^{t}\left[  \sup\limits_{s\in\left[  t,T\right]  }\left\vert
z^{\varepsilon,v}\left(  s\right)  \right\vert ^{2k}\right]  =O\left(
\varepsilon^{2k}\right)  ,\tag{3.16}\\
&  \mathbb{E}^{t}\left[  \sup\limits_{s\in\left[  t,T\right]  }\left\vert
y^{\varepsilon,v}\left(  s\right)  +z^{\varepsilon,v}\left(  s\right)
\right\vert ^{2k}\right]  =O\left(  \varepsilon^{k}\right)  . \tag{3.17}%
\end{align}

In addition, we have the following equality%

\begin{align}
&  J\left(  t,\hat{X}\left(  t\right)  ,u^{\varepsilon}\left(  \cdot\right)
\right)  -J\left(  t,\hat{X}\left(  t\right)  ,\hat{u}\left(  \cdot\right)
\right) \nonumber\\
&  ={\displaystyle\int_{t}^{t+\varepsilon}} \mathbb{E}^{t}\left[  \left\langle
\mathcal{H}\left(  s;t\right)  ,v\right\rangle +\frac{1}{2}\left\langle
\mathcal{A}\left(  s;t\right)  v,v\right\rangle \right]  ds +o\left(
\varepsilon\right)  . \tag{3.18}%
\end{align}

\end{proposition}

\begin{proof}
See the Appendix.
\end{proof}

Now, we present the following technical lemma needed later. The proof follows
an argument adapted from Hu el al. \cite{Huetal2015},

\begin{lemma}
Under assumptions \textbf{(H1)-(H3)}, the following two statements are equivalent

\begin{itemize}
\item[i)] $\underset{\varepsilon\downarrow0}{\lim}\dfrac{1}{\varepsilon
}{\displaystyle\int_{t}^{t+\varepsilon}} \mathbb{E}^{t}\left[  \mathcal{H}%
\left(  s;t\right)  \right]  ds=0,$ $d\mathbb{P}-a.s,$ $\forall t\in\left[
0,T\right]  .$

\item[ii)] $\mathcal{H}\left(  t;t\right)  =0,$ $d\mathbb{P}-a.s,$ $dt-a.e.$
\end{itemize}
\end{lemma}

\begin{proof}
See the Appendix.
\end{proof}

\noindent\textbf{Proof of Theorem 3.2.} Given an admissible strategy%

\[
\hat{u}\left(  \cdot\right)  \in\mathcal{L}_{\mathcal{F}}^{1}\left(
0,T;\mathbb{R}\right)  \times\mathcal{M}_{\mathcal{F}}^{2}\left(
0,T;\mathbb{R}^{d}\right)  ,
\]
for which $\left(  3.11\right)  $ holds$,$ according to Lemma 3.4, we have,
for any $t\in\left[  0,T\right]  ,$%

\[
\lim\limits_{\varepsilon\downarrow0}\dfrac{1}{\varepsilon} {\displaystyle\int
_{t}^{t+\varepsilon}}\mathbb{E}^{t}\left[  \mathcal{H}\left(  s;t\right)
\right]  ds=0, \, a.s.
\]

Then, from $\left(  3.18\right)  ,$ for any $t\in\left[  0,T\right]  $ and for
any $\mathbb{R}^{d+1}-$valued, $\mathcal{F}_{t}-$measurable and bounded random
variable $v,$%

\begin{align*}
&  \lim_{\varepsilon\downarrow0}\frac{1}{\varepsilon}\left\{  J\left(
t,\hat{X}\left(  t\right)  ,u^{\varepsilon}\left(  \cdot\right)  \right)
-J\left(  t,\hat{X}\left(  t\right)  ,\hat{u}\left(  \cdot\right)  \right)
\right\} \\
&  =\lim_{\varepsilon\downarrow0}\frac{1}{\varepsilon}{\displaystyle \int
_{t}^{t+\varepsilon}} \left\{  \left\langle \mathbb{E}^{t}\left[
\mathcal{H}\left(  s;t\right)  \right]  ,v\right\rangle ds+\frac{1}%
{2}\left\langle \mathbb{E}^{t}\left[  \mathcal{A}\left(  s;t\right)  \right]
v,v\right\rangle \right\}  ds,\\
&  =\frac{1}{2}\lim_{\varepsilon\downarrow0}\frac{1}{\varepsilon
}{\displaystyle\int_{t}^{t+\varepsilon}} \left\langle \mathbb{E}^{t}\left[
\mathcal{A}\left(  s;t\right)  \right]  v,v\right\rangle ds,\\
&  \leq0,
\end{align*}
where we have used in the last inequality the fact that, under the concavity
condition of $\varphi\left(  \cdot\right)  $ and $h\left(  \cdot\right)
$,\textbf{\ }it follows $\left\langle \mathcal{A}\left(  s;t\right)
v,v\right\rangle \leq0.$ Hence $\hat{u}\left(  \cdot\right)  $ is an
equilibrium strategy.

Conversely, assume that $\hat{u}\left(  \cdot\right)  $ is an equilibrium
strategy. Then, by $\left(  3.2\right)  $ together with $\left(  3.18\right)
,$ for any $\left(  t,u\right)  \in\left[  0,T\right]  \times\mathbb{R}%
^{d+1},$ the following inequality holds:
\begin{equation}
\lim_{\varepsilon\downarrow0}\left\langle \frac{1}{\varepsilon}%
{\displaystyle\int_{t}^{t+\varepsilon}} \mathbb{E}^{t}\left[  \mathcal{H}%
\left(  s;t\right)  \right]  ds,u\right\rangle +\frac{1}{2}\lim_{\varepsilon
\downarrow0}\left\langle \frac{1}{\varepsilon}{\displaystyle\int
_{t}^{t+\varepsilon}} \mathbb{E}^{t}\left[  \mathcal{A}\left(  s;t\right)
\right]  ds \, u,u\right\rangle \leq0. \tag{3.19}%
\end{equation}

Now, we define $\forall\left(  t,u\right)  \in\left[  0,T\right]
\times\mathbb{R}^{d+1},$%

\[
\Phi\left(  t,u\right)  =\lim_{\varepsilon\downarrow0}\frac{1}{\varepsilon
}\left\langle {\displaystyle\int_{t}^{t+\varepsilon}} \mathbb{E}^{t}\left[
\mathcal{H}\left(  s;t\right)  \right]  ds,u\right\rangle +\frac{1}{2}%
\lim_{\varepsilon\downarrow0}\left\langle \frac{1}{\varepsilon}%
{\displaystyle\int_{t}^{t+\varepsilon}} \mathbb{E}^{t}\left[  \mathcal{A}%
\left(  s;t\right)  \right]  ds \, u,u\right\rangle .
\]

Clearly $\Phi\left(  \cdot, \cdot\right)  $ is well defined. In fact, it is a
second order polynomial in terms of the components of vector $u.$ Easy
manipulations show that the inequality $\left(  3.19\right)  $ is equivalent to%

\begin{equation}
\Phi\left(  t,0\right)  =\max_{u\in\mathbb{R}^{d+1}}\Phi\left(  t,u\right)  ,
\text{ }d\mathbb{P}-a.s,\forall t\in\left[  0,T\right]  . \tag{3.20}%
\end{equation}
And it is easy to see that the maximum condition $\left(  3.20\right)  $ leads
to the following condition: $\forall t\in\left[  0,T\right]  , $%

\begin{equation}
\Phi_{u}\left(  t,0\right)  =\lim_{\varepsilon\downarrow0}\frac{1}%
{\varepsilon}{\displaystyle\int_{t}^{t+\varepsilon}} \mathbb{E}^{t}\left[
\mathcal{H}\left(  s;t\right)  \right]  ds=0,\text{ }d\mathbb{P}-a.s.
\tag{3.21}%
\end{equation}

According to Lemma 3.4, the expression $\left(  3.11\right)  $ follows
immediately. This completes the proof.\eop

\subsection{A characterization of equilibrium strategies by verification
argument}

In classical (time-consistent) stochastic control theory the sufficient
condition of optimality is of significant importance for computing optimal
controls. It says that if an admissible control satisfies the maximum
condition of the Hamiltonian function, then the control is indeed optimal for
the stochastic control problem. This allows one to solve examples of optimal
control problems where one can find, a smooth solution to the associated
adjoint equation.

It is worth mentioning also that, the assumption \textbf{(H3)} is too strong
to apply the theorem 3.2. to some important problems in the practice. For
example, the power utility function does not satisfy this assumption. The aim
of the following theorem is to characterize the open-loop equilibrium pair by
a sufficient condition of equilibrium. In order to overcome the technical
difficulties mentioned by the hypothesis \textbf{(H3)}, let us introduce
further conditions about the utility functions.

\begin{enumerate}
\item[\textbf{(H4)}] The maps $\varphi\left(  \cdot\right)  ,h\left(
\cdot\right)  $ are continuously differentiable and the first order
derivatives $\varphi_{x}\left(  \cdot\right)  ,h_{x}\left(  \cdot\right)  $
are continuous.
\end{enumerate}

\begin{theorem}
Let \textbf{(H1), (H2) }and \textbf{(H4)} hold. Given an admissible strategy
$\hat{u}\left(  \cdot\right)  \in\mathcal{L}_{\mathcal{F}}^{1}\left(
0,T;\mathbb{R}\right)  \times\mathcal{M}_{\mathcal{F}}^{2}\left(
0,T;\mathbb{R}^{d}\right)  $, let for any $t\in\left[  0,T\right]  ,$ the
processus
\[
\left(  p\left(  \cdot;t\right)  ,q\left(  \cdot;t\right)  \right)
\in\mathcal{L}_{\mathcal{F}}^{2}\left(  t,T;\mathbb{R}\right)  \times
\mathcal{M}_{\mathcal{F}}^{2}\left(  t,T;\mathbb{R}^{d}\right)
\]
be the unique solution to the BSDE $\left(  3.4\right)  $. Then, $\hat
{u}\left(  \cdot\right)  $ is an equilibrium consumption-investment strategy,
if the following condition holds%
\begin{equation}
\mathcal{H}\left(  t;t\right)  =0,\text{ }d\mathbb{P-}a.s.,\text{ }dt-a.e.
\tag{3.22}%
\end{equation}

\end{theorem}

\textbf{Proof. }Suppose that $\hat{u}(\cdot)$ is an admissible control for
which the condition (3.22) holds, in addition for any $t\in\lbrack0,T]$ and
$\varepsilon\in\lbrack0,T-t)$, we consider $u^{\varepsilon}\left(
\cdot\right)  $ by $\left(  3.1\right)  ,$ then we have the following
difference%
\begin{align*}
&  J\left(  t,\hat{X}\left(  t\right)  ,\hat{u}\left(  \cdot\right)  \right)
-J\left(  t,\hat{X}\left(  t\right)  ,u^{\varepsilon}\left(  \cdot\right)
\right) \\
&  =\mathbb{E}^{t}\left[  {\displaystyle\int_{t}^{T}}\lambda\left(
s-t\right)  \left(  \varphi\left(  \Gamma^{\top}\hat{u}\left(  s\right)
\right)  -\varphi\left(  \Gamma^{\top}u^{\varepsilon}\left(  s\right)
\right)  \right)  ds+\lambda\left(  T-t\right)  \left(  h\left(  \hat
{X}\left(  T\right)  \right)  -h\left(  \hat{X}^{\varepsilon}\left(  T\right)
\right)  \right)  \right]  .
\end{align*}

Noting that, by the concavity of $h\left(  \cdot\right)  $, we have%
\[%
\begin{array}
[c]{ll}%
\mathbb{E}^{t}\left[  \lambda\left(  T-t\right)  \left(  h\left(  \hat
{X}\left(  T\right)  \right)  -h\left(  \hat{X}^{\varepsilon}\left(  T\right)
\right)  \right)  \right]  & \geq\mathbb{E}^{t}\left[  \lambda\left(
T-t\right)  \left(  \hat{X}\left(  T\right)  -\hat{X}^{\varepsilon}\left(
T\right)  \right)  ^{T}h_{x}\left(  \hat{X}\left(  T\right)  \right)  \right]
,
\end{array}
\]

Accordingly, by the terminal condition in the BSDE (3.4) we obtain that%
\begin{align}
&  J\left(  t,\hat{X}\left(  t\right)  ,\hat{u}\left(  \cdot\right)  \right)
-J\left(  t,\hat{X}\left(  t\right)  ,u^{\varepsilon}\left(  \cdot\right)
\right) \nonumber\\
&  \geq\mathbb{E}^{t}\left[  {\displaystyle\int_{t}^{T}}\lambda\left(
s-t\right)  \left(  \varphi\left(  \Gamma^{\top}\hat{u}\left(  s\right)
\right)  -\varphi\left(  \Gamma^{\top}u^{\varepsilon}\left(  s\right)
\right)  \right)  ds+\left(  \hat{X}\left(  T\right)  -\hat{X}^{\varepsilon
}\left(  T\right)  \right)  ^{T}p\left(  T;t\right)  \right]  . \tag{3.23}%
\end{align}

By applying Ito's formula to $s\mapsto\left(  \hat{X}\left(  s\right)
-\hat{X}^{\varepsilon}\left(  s\right)  \right)  ^{T}p\left(  s;t\right)  $ on
$\left[  t,T\right]  $, we get%
\begin{equation}
\mathbb{E}^{t}\left[  \left(  \hat{X}\left(  T\right)  -\hat{X}^{\varepsilon
}\left(  T\right)  \right)  ^{T}p\left(  T;t\right)  \right]  =\mathbb{E}%
^{t}\left[  {\displaystyle\int_{t}^{T}}\left(  \hat{u}\left(  s\right)
-u^{\varepsilon}\left(  s\right)  \right)  ^{T}\left(  B\left(  s\right)
p\left(  s;t\right)  +D\left(  s\right)  \widetilde{q}\left(  s;t\right)
\right)  ds\right]  . \tag{3.24}%
\end{equation}

By the concavity of $\varphi\left(  \cdot\right)  $, we find%
\begin{align}
&  \mathbb{E}^{t}\left[  {\displaystyle\int_{t}^{T}}\lambda\left(  s-t\right)
\left(  \varphi\left(  \Gamma^{\top}\hat{u}\left(  s\right)  \right)
-\varphi\left(  \Gamma^{\top}u^{\varepsilon}\left(  s\right)  \right)
\right)  ds\right] \nonumber\\
&  \geq\mathbb{E}^{t}\left[  {\displaystyle\int_{t}^{T}}\lambda\left(
s-t\right)  \left\langle \varphi_{c}\left(  \Gamma^{\top}\hat{u}\left(
s\right)  \right)  \Gamma,\hat{u}\left(  s\right)  -u^{\varepsilon}\left(
s\right)  \right\rangle ds\right]  , \tag{3.25}%
\end{align}

By taking $\left(  3.24\right)  $ and $\left(  3.25\right)  $ in $\left(
3.23\right)  ,$ it follows that%
\begin{align*}
&  J\left(  t,\hat{X}\left(  t\right)  ,u^{\varepsilon}\left(  \cdot\right)
\right)  -J\left(  t,\hat{X}\left(  t\right)  ,\hat{u}\left(  \cdot\right)
\right) \\
&  \leq\mathbb{E}^{t}\left[  {\displaystyle\int_{t}^{T}}\left\langle B\left(
s\right)  p\left(  s;t\right)  +D\left(  s\right)  \widetilde{q}\left(
s;t\right)  +\lambda\left(  s-t\right)  \varphi_{c}\left(  \Gamma^{\top}%
\hat{u}\left(  s\right)  \right)  \Gamma,u^{\varepsilon}\left(  s\right)
-\hat{u}\left(  s\right)  \right\rangle ds\right] \\
&  =\mathbb{E}^{t}\left[  {\displaystyle\int_{t}^{t+\varepsilon}}\left\langle
\mathcal{H}\left(  s;t\right)  ,v\right\rangle ds\right]  .
\end{align*}
Now dividing both sides by $\varepsilon$ and taking the limit when
$\varepsilon$ vanishes, by Lemma 3.4, we conclude tha $\hat{u}($%
\textperiodcentered$)$ is an equilibrium control.\eop

\begin{remark}
The purpose of the sufficient condition of optimality is to find an optimal
control by computing the difference $J\left(  \hat{u}\left(  \cdot\right)
\right)  -J\left(  u\left(  \cdot\right)  \right)  $ in terms of the
Hamiltonian function,\ where $u\left(  \cdot\right)  $ be an arbitrary
admissible control.\ Here, the spike variation perturbation $\left(
3.1\right)  $\ plays a key role in deriving the sufficient condition for
equilibrium strategies, which reduces to the computation of the difference
$J\left(  t,\hat{X}\left(  t\right)  ,\hat{u}\left(  \cdot\right)  \right)
-J\left(  t,\hat{X}\left(  t\right)  ,u^{\varepsilon}\left(  \cdot\right)
\right)  $, without the necessity to achieving the second order expansion in
the spike variation.
\end{remark}

\section{Equilibrium when the coefficients are deterministic}

Theorem 3.5. shows that one can obtain equilibrium consumption-investment
strategies by solving a system of FBSDEs which is not standard since the
\textquotedblleft flow\textquotedblright\ of the unknown process $\left(
p\left(  \cdot;t\right)  ,q\left(  \cdot;t\right)  \right)  _{t\in\left[
0,T\right]  }$ is involved. Moreover, there is an additional constraint that
act on the \textquotedblleft diagonal\textquotedblright\ (i.e. when $s=t$) of
the flow. As far as we know, the explicitly solvability of this type of
equations remains an open problem, except for some particular form of the
utility functions. However, we are able to solve quite thoroughly this problem
when the parameters $r\left(  \cdot\right)  $ and $\sigma\left(  \cdot\right)
$ are deterministic functions. In this section we define what we mean by an
equilibrium rule, and then we derive a parabolic backward PDE. Our PDE is
comparable with the one obtained in \cite{Solano} and \cite{EkePir}$,$ for
some particular discount functions in finite horizon with different utility functions.

In this section, let us look at the Merton's portfolio problem with general
discounting and deterministic parameters. At first, we consider the following
parabolic backward partial differential equation
\begin{equation}
\left\{
\begin{array}
[c]{l}%
\theta_{t}\left(  t,x\right)  +\theta_{x}\left(  t,x\right)  \left(
r_{0}\left(  t\right)  x-r\left(  t\right)  ^{\top}\Sigma\left(  t\right)
r\left(  t\right)  \dfrac{\theta\left(  t,x\right)  }{\theta_{x}\left(
t,x\right)  }-\mathcal{I}\left(  \lambda\left(  T-t\right)  \theta\left(
t,x\right)  \right)  \right) \\
+\dfrac{1}{2}\theta_{xx}\left(  t,x\right)  r\left(  t\right)  ^{\top}%
\Sigma\left(  t\right)  r\left(  t\right)  \left(  \dfrac{\theta\left(
t,x\right)  }{\theta_{x}\left(  t,x\right)  }\right)  ^{2}+\theta\left(
t,x\right)  r_{0}\left(  t\right)  =0,\text{ }\left(  t,x\right)  \in\left[
0,T\right]  \times\mathbb{R},\\
\theta\left(  T,x\right)  =h_{x}\left(  x\right)  ,
\end{array}
\right.  \tag{4.1}%
\end{equation}
where we denote by $\mathcal{I}\left(  \cdot\right)  $ the inverse function of
the strictly decreasing marginal derivative utility $\varphi_{c}\left(
\cdot\right)  $ and $\Sigma\left(  s\right)  \equiv\left(  \sigma\left(
s\right)  \sigma\left(  s\right)  ^{\top}\right)  ^{-1}$.

We have the following verification theorem

\begin{theorem}
Let \textbf{(H1), (H2) }and\textbf{ (H4) }hold. If there exists a classical
solution
\[
\theta\left(  \cdot,\cdot\right)  \in\mathcal{C}^{1,2}\left(  (0,T)\times
\mathbb{R},\mathbb{R}\right)  \cap\mathcal{C}\left(  \left[  0,T\right]
\times\mathbb{R},\mathbb{R}\right)
\]
of the PDE $\left(  4.1\right)  $ such that the stochastic differential
equation%
\begin{equation}
\left\{
\begin{array}
[c]{l}%
d\hat{X}\left(  s\right)  =\left\{  r_{0}\left(  s\right)  \hat{X}\left(
s\right)  -r\left(  s\right)  ^{\top}\Sigma\left(  s\right)  r\left(
s\right)  \dfrac{\theta\left(  s,\hat{X}\left(  s\right)  \right)  }%
{\theta_{x}\left(  s,\hat{X}\left(  s\right)  \right)  }-\mathcal{I}\left(
\lambda\left(  T-s\right)  \theta\left(  s,\hat{X}\left(  s\right)  \right)
\right)  \right\}  ds\\
\text{ \ \ \ \ \ \ \ \ \ \ \ \ \ \ \ \ \ }-\dfrac{\theta\left(  s,\hat
{X}\left(  s\right)  \right)  }{\theta_{x}\left(  s,\hat{X}\left(  s\right)
\right)  }r\left(  s\right)  ^{\top}\Sigma\left(  s\right)  \sigma\left(
s\right)  dW\left(  s\right)  ,\text{ }s\in\left[  0,T\right]  ,\\
\hat{X}\left(  0\right)  =x_{0},
\end{array}
\right.  \tag{4.2}%
\end{equation}
has a unique solution $\hat{X}\left(  \cdot\right)  ,$ in which the following
estimate holds
\[
\mathbb{E}\left[  \underset{0\leq t\leq T}{\sup}\left\vert X\left(  t\right)
\right\vert ^{2}\right]  \leq K\left(  1+\left\vert x_{0}\right\vert
^{2}\right)  .
\]

Then, the equilibrium consumption-investment strategy $\hat{u}\left(
\cdot\right)  =\left(  \hat{c}\left(  \cdot\right)  ,\hat{u}_{I}\left(
\cdot\right)  ^{\top}\right)  ^{\top}$ is given by%
\begin{align}
&  \hat{c}\left(  t\right)  =\mathcal{I}\left(  \lambda\left(  T-t\right)
\theta\left(  t,\hat{X}\left(  t\right)  \right)  \right)  ,\text{
}dt-a.e.,\tag{4.3}\\
&  \hat{u}_{I}\left(  t\right)  =-\Sigma\left(  t\right)  r\left(  t\right)
\frac{\theta\left(  t,\hat{X}\left(  t\right)  \right)  }{\theta_{x}\left(
t,\hat{X}\left(  t\right)  \right)  },\text{ }dt-a.e. \tag{4.4}%
\end{align}

\end{theorem}

\textbf{Proof.} Suppose that $\hat{u}\left(  \cdot\right)  =\left(  \hat
{c}\left(  \cdot\right)  ,\hat{u}_{I}\left(  \cdot\right)  ^{\top}\right)
^{\top}$ is an equilibrium control and denote by $\hat{X}\left(  \cdot\right)
$ the corresponding wealth process. Then, in view of Theorem $3.5$, there
exist an adapted process $\left(  \hat{X}\left(  \cdot\right)  ,\left(
p\left(  \cdot;t\right)  ,q\left(  \cdot;t\right)  \right)  _{t\in\left[
0,T\right]  }\right)  ,$ solution of the following flow of forward-backward
SDEs, parametrized by $t\in\left[  0,T\right]  :$%
\begin{equation}
\left\{
\begin{array}
[c]{l}%
dX\left(  s\right)  =\left\{  r_{0}\left(  s\right)  \hat{X}\left(  s\right)
+\hat{u}_{I}\left(  s\right)  ^{\top}r\left(  s\right)  -\hat{c}\left(
s\right)  \right\}  ds+\hat{u}_{I}\left(  s\right)  ^{\top}\sigma\left(
s\right)  dW\left(  s\right)  ,\text{ }s\in\left[  t,T\right]  ,\\
dp\left(  s;t\right)  =-r_{0}\left(  s\right)  p\left(  s;t\right)
ds+q\left(  s,t\right)  ^{\top}dW\left(  s\right)  ,\text{ }0\leq t\leq s\leq
T,\\
\hat{X}\left(  0\right)  =x_{0},\text{ }p\left(  T;t\right)  =\lambda\left(
T-t\right)  h_{x}\left(  \hat{X}\left(  T\right)  \right)  ,\text{ }%
t\in\left[  0,T\right]  ,
\end{array}
\right.  \tag{4.5}%
\end{equation}
with conditions%
\begin{align}
&  -p\left(  t;t\right)  +\varphi_{c}\left(  \hat{c}\left(  t\right)  \right)
=0,\text{ }dt-a.e.,\tag{4.6}\\
&  p\left(  t;t\right)  r\left(  t\right)  +\sigma\left(  t\right)  q\left(
t;t\right)  =0,\text{ }dt-a.e. \tag{4.7}%
\end{align}

From the terminal condition in the first order adjoint process we consider the
following Ansatz%
\begin{equation}
p\left(  s;t\right)  =\lambda\left(  T-t\right)  \mathcal{V}\left(  s,\hat
{X}\left(  s\right)  \right)  ,\text{ }\forall\text{ }0\leq t\leq s\leq T,
\tag{4.8}%
\end{equation}
for some deterministic function $\mathcal{V}\left(  \cdot, \cdot\right)
\in\mathcal{C}^{1,2}\left(  \left[  0,T\right]  \times\mathbb{R},
\mathbb{R}\right)  $ such that $\mathcal{V}\left(  T, \cdot\right)
=h_{x}\left(  \cdot\right)  .$

Applying It\^{o}'s formula to $\left(  4.8\right)  $, it yields%
\begin{align}
dp\left(  s;t\right)   &  =\lambda\left(  T-t\right)  \left\{  \mathcal{V}%
_{s}\left(  s,\hat{X}\left(  s\right)  \right)  +\mathcal{V}_{x}\left(
s,\hat{X}\left(  s\right)  \right)  \left(  \hat{X}\left(  s\right)
r_{0}\left(  s\right)  +\hat{u}_{I}\left(  s\right)  ^{\top}r\left(  s\right)
-\hat{c}\left(  s\right)  \right)  \right. \nonumber\\
&  \text{ \ \ \ \ \ \ \ \ \ \ \ \ \ \ \ \ \ \ \ \ \ \ \ \ }\left.  +\frac
{1}{2}\mathcal{V}_{xx}\left(  s,\hat{X}\left(  s\right)  \right)  \hat{u}%
_{I}\left(  s\right)  ^{\top}\sigma\left(  s\right)  \sigma\left(  s\right)
^{\top}\hat{u}_{I}\left(  s\right)  \right\}  ds\nonumber\\
&  \text{ \ \ \ \ \ \ }+\lambda\left(  T-t\right)  \mathcal{V}_{x}\left(
s,\hat{X}\left(  s\right)  \right)  \hat{u}_{I}\left(  s\right)  ^{\top}%
\sigma\left(  s\right)  dW\left(  s\right)  . \tag{4.9}%
\end{align}

Next, comparing the $ds$ term in $\left(  4.9\right)  $ by the ones in the
second equation in $\left(  4.5\right)  ,$ we deduce that%
\begin{align}
&  \mathcal{V}_{s}\left(  s,\hat{X}\left(  s\right)  \right)  +\mathcal{V}%
_{x}\left(  s,\hat{X}\left(  s\right)  \right)  \left(  \hat{X}\left(
s\right)  r_{0}\left(  s\right)  +\hat{u}_{I}\left(  s\right)  ^{\top}r\left(
s\right)  -\hat{c}\left(  s\right)  \right) \nonumber\\
&  +\frac{1}{2}\mathcal{V}_{xx}\left(  s,\hat{X}\left(  s\right)  \right)
\hat{u}_{I}\left(  s\right)  ^{\top}\sigma\left(  s\right)  \sigma\left(
s\right)  ^{\top}\hat{u}_{I}\left(  s\right)  =-r_{0}\left(  s\right)
\mathcal{V}\left(  s,\hat{X}\left(  s\right)  \right)  , \tag{4.10}%
\end{align}
and by comparing the $dW\left(  s\right)  $ terms we also get
\begin{equation}
q\left(  s,t\right)  =\lambda\left(  T-t\right)  \mathcal{V}_{x}\left(
s,\hat{X}\left(  s\right)  \right)  \sigma\left(  s\right)  ^{\top}\hat{u}%
_{I}\left(  s\right)  . \tag{4.11}%
\end{equation}

We put the above expressions of $p\left(  s;t\right)  $ and $q\left(
s;t\right)  $ at $s=t$\ into $\left(  4.6\right)  $ and $\left(  4.7\right)
,$ then%
\begin{equation}
\lambda\left(  T-t\right)  \mathcal{V}\left(  t,\hat{X}\left(  t\right)
\right)  -\varphi_{c}\left(  \hat{c}\left(  t\right)  \right)  =0, \tag{4.12}%
\end{equation}
and%
\begin{equation}
\mathcal{V}_{x}\left(  t,\hat{X}\left(  t\right)  \right)  \sigma\left(
t\right)  \sigma\left(  t\right)  ^{\top}\hat{u}_{I}\left(  t\right)
=-r\left(  t\right)  \mathcal{V}\left(  t,\hat{X}\left(  t\right)  \right)  ,
\tag{4.13}%
\end{equation}
which leads to the following representation%
\begin{align}
&  \hat{c}\left(  t\right)  =\mathcal{I}\left(  \lambda\left(  T-t\right)
\mathcal{V}\left(  t,\hat{X}\left(  t\right)  \right)  \right)  ,\text{
}dt-a.e.,\tag{4.14}\\
&  \hat{u}_{I}\left(  t\right)  =-\Sigma\left(  t\right)  r\left(  t\right)
\frac{\mathcal{V}\left(  t,\hat{X}\left(  t\right)  \right)  }{\mathcal{V}%
_{x}\left(  t,\hat{X}\left(  t\right)  \right)  },\text{ }dt-a.e. \tag{4.15}%
\end{align}

Then by taking expressions $\left(  4.14\right)  $ and $\left(  4.15\right)  $
into $\left(  4.10\right)  $, this suggests that $\mathcal{V}\left(  \cdot,
\cdot\right)  $ coincides with the solution of the PDE $\left(  4.1\right)  $,
evaluated along the trajectory $\hat{X}\left(  t\right)  ,$ solution of the
state equation$.$

\begin{remark}
Equation $\left(  4.1\right)  $ is comparable with the one in
Mar\'{\i}n-Solano \& Navas \cite{Solano} and Ekland \& Pirvu \cite{EkePir}, in
which the equilibrium is defined within the class of feedback controls.
\end{remark}

\begin{remark}
Theorem 4.1. enables us to derive a suitable equilibrium strategie $\hat
{u}_{I}\left(  t\right)  $ as well as $\hat{c}\left(  t\right)  $, at each
$t\in\lbrack0,T]$, this permits us to derive directly an explicit expression
of equilibrium control in the cases of power, logarithmic and exponential
utility functions. While the duality approach \cite{Yus} permits to
characterize\ a stochastic equilibrium solution in terms of a complicated
FBSDE system of a closed form; it does not provide an explicit represntation.
\end{remark}

\section{Special utility functions}

Equilibrium investment-consumption strategies for Merton's portfolio problem
with general discounting and deterministic parameters have been studied in
\cite{Solano}, \cite{EkePir} and \cite{Yong2012}, among others, in different
frameworks. In this section, we discuss some special cases in which the
function $\theta\left(  \cdot, \cdot\right)  $ may be separated into functions
of time and state variables. Then, one needs only to solve a system of ODEs in
order to completely determine the equilibrium strategies. We will compare our
results with some existing ones in the literature.

\subsection{Power utility function}

To make the problem $\left(  2.8\right)  -\left(  2.9\right)  $ explicitly
solvable, we consider power utility functions for the running and terminal
costs. That is, $\varphi\left(  c\right)  =\frac{c^{\gamma}}{\gamma}$ and
$h\left(  x\right)  =a\frac{x^{\gamma}}{\gamma},$ with $a>0$ and $\gamma
\in\left(  0,1\right)  .$ In this case the PDE $\left(  4.1\right)  $ reduces
to%
\[
\left\{
\begin{array}
[c]{l}%
\theta_{t}\left(  t,x\right)  +\theta_{x}\left(  t,x\right)  \left(
r_{0}\left(  t\right)  x-r\left(  t\right)  ^{\top}\Sigma\left(  t\right)
r\left(  t\right)  \dfrac{\theta\left(  t,x\right)  }{\theta_{x}\left(
t,x\right)  }-\dfrac{\lambda\left(  T-t\right)  ^{1-\gamma}}{\theta\left(
t,x\right)  ^{\gamma-1}}\right) \\
+\dfrac{1}{2}\theta_{xx}\left(  t,x\right)  r\left(  t\right)  ^{\top}%
\Sigma\left(  t\right)  r\left(  t\right)  \left(  \dfrac{\theta\left(
t,x\right)  }{\theta_{x}\left(  t,x\right)  }\right)  ^{2}+r_{0}\left(
t\right)  \theta\left(  t,x\right)  =0,\text{ }\left(  t,x\right)  \in\left[
0,T\right]  \times\mathbb{R},\\
\theta\left(  T,x\right)  =ax^{\gamma-1}.
\end{array}
\right.
\]

From the terminal condition, we consider the following trial solution%
\[
\theta\left(  s,x\right)  =a\Pi\left(  s\right)  x^{\gamma-1},
\]
for some deterministic function $\Pi\left(  \cdot\right)  \in C^{1}\left(
\left[  0,T\right]  ,\mathbb{R}\right)  $ with the terminal condition
$\Pi\left(  T\right)  =1.$ Then by substituting in $\left(  4.1\right)  ,$ we
obtain%
\begin{equation}
\left\{
\begin{array}
[c]{l}%
\Pi_{t}\left(  t\right)  +\left(  K\left(  t\right)  +Q\left(  t\right)
\Pi\left(  t\right)  ^{\frac{1}{\gamma-1}}\right)  \Pi\left(  t\right)
=0,\text{ for }t\in\left[  0,T\right]  ,\\
\Pi\left(  T\right)  =1.
\end{array}
\right.  \tag{5.1}%
\end{equation}
where%
\begin{equation}
K\left(  t\right)  \equiv\gamma r_{0}\left(  t\right)  +\frac{1}{2}%
\dfrac{\gamma}{\left(  1-\gamma\right)  }r\left(  t\right)  ^{\top}%
\Sigma\left(  t\right)  r\left(  t\right)  \text{,} \tag{5.2}%
\end{equation}
and%
\begin{equation}
Q\left(  t\right)  \equiv\left(  1-\gamma\right)  \left(  a\lambda\left(
T-t\right)  \right)  ^{\frac{1}{\gamma-1}}. \tag{5.3}%
\end{equation}

It remains to determine the function $\Pi\left(  \cdot\right)  .$ First, by
the change of variable
\begin{equation}
\Pi\left(  t\right)  =y\left(  t\right)  ^{\left(  1-\gamma\right)  },\text{
for }t\in\left[  0,T\right]  , \tag{5.4}%
\end{equation}
we find that $y\left(  \cdot\right)  $ should solve the following ODE%
\[
\left\{
\begin{array}
[c]{l}%
y_{t}\left(  t\right)  -\dfrac{K\left(  t\right)  }{\left(  \gamma-1\right)
}y\left(  t\right)  -\dfrac{Q\left(  t\right)  }{\left(  \gamma-1\right)
}=0,\text{ for }t\in\left[  0,T\right]  ,\\
y\left(  T\right)  =1.
\end{array}
\right.
\]

A variation of constant formula yields to%
\[
y\left(  t\right)  =\left(  1-\int_{t}^{T}\dfrac{Q\left(  \tau\right)
}{\left(  \gamma-1\right)  }e^{{\displaystyle\int_{\tau}^{T}}\dfrac{K\left(
l\right)  }{\left(  \gamma-1\right)  }dl}d\tau\right)  \exp\left(
-{\displaystyle\int_{t}^{T}}\dfrac{K\left(  \tau\right)  }{\left(
\gamma-1\right)  }d\tau\right)  ,\text{ for }t\in\left[  0,T\right]  ,
\]
and subsequently we obtain%

\[
\Pi\left(  t\right)  =\left(  1-\int_{t}^{T}\dfrac{Q\left(  \tau\right)
}{\left(  \gamma-1\right)  }e^{{\displaystyle\int_{\tau}^{T}} \dfrac{K\left(
l\right)  }{\left(  \gamma-1\right)  }dl}d\tau\right)  ^{1-\gamma}\exp\left(
{\displaystyle\int_{t}^{T}} K\left(  \tau\right)  d\tau\right)  ,\text{ for
}t\in\left[  0,T\right]  .
\]

In view of Theorem 4.1, the representation of the\ Nash
equilibrium\ strategies $\left(  4.3\right)  $-$\left(  4.4\right)  $ gives
\begin{align}
&  \hat{c}\left(  t\right)  =\left(  a\lambda\left(  T-t\right)  \Pi\left(
t\right)  \right)  ^{\frac{1}{\gamma-1}}\hat{X}\left(  t\right)  ,\text{
}dt-a.e.,\tag{5.5}\\
&  \hat{u}_{I}\left(  t\right)  =\Sigma\left(  t\right)  r\left(  t\right)
\frac{\hat{X}\left(  t\right)  }{\left(  1-\gamma\right)  },\text{ }dt-a.e.
\tag{5.6}%
\end{align}

This consumption--investment strategy determines a wealth process given by%
\[%
\begin{array}
[c]{cl}%
X\left(  t\right)  = & x_{0}+%
{\displaystyle\int_{0}^{t}}
\left\{  r_{0}\left(  s\right)  +\dfrac{1}{\left(  1-\gamma\right)  }r\left(
s\right)  ^{\top}\Sigma\left(  s\right)  r\left(  s\right)  -\left(
a\lambda\left(  T-s\right)  \Pi\left(  s\right)  \right)  ^{\frac{1}{\gamma
-1}}\right\}  \hat{X}\left(  s\right)  ds\\
& \text{ \ \ \ }+%
{\displaystyle\int_{0}^{t}}
\dfrac{\hat{X}\left(  s\right)  }{\left(  1-\gamma\right)  }r\left(  s\right)
^{\top}\Sigma\left(  s\right)  \sigma\left(  s\right)  dW\left(  s\right)
,\text{ }t\in\left[  0,T\right]  ,
\end{array}
\]

The abouve solution is comparable with the one obtained by Mar\'{\i}n-Solano
\& Navas \cite{Solano}, Ekland \& Pirvu \cite{EkePir} and Yong \cite{Yong2012}.

\subsection{Logarithmic utility function}

Now, let us analyze the case where $\varphi\left(  c\right)  =\ln\left(
c\right)  $ and $h\left(  x\right)  =a\ln\left(  x\right)  ,$ with $a>0.$ In
this case, the PDE $\left(  4.1\right)  $ reduces to%
\begin{equation}
\left\{
\begin{array}
[c]{l}%
\theta_{t}\left(  t,x\right)  +\theta_{x}\left(  t,x\right)  \left(
r_{0}\left(  t\right)  x-r\left(  t\right)  ^{\top}\Sigma\left(  t\right)
r\left(  t\right)  \dfrac{\theta\left(  t,x\right)  }{\theta_{x}\left(
t,x\right)  }-\left(  \lambda\left(  T-t\right)  \theta\left(  t,x\right)
\right)  ^{-1}\right) \\
+\dfrac{1}{2}\theta_{xx}\left(  t,x\right)  r\left(  t\right)  ^{\top}%
\Sigma\left(  t\right)  r\left(  t\right)  \left(  \dfrac{\theta\left(
t,x\right)  }{\theta_{x}\left(  t,x\right)  }\right)  ^{2}+r_{0}\left(
t\right)  \theta\left(  t,x\right)  =0,\text{ }\left(  t,x\right)  \in\left[
0,T\right]  \times\mathbb{R},\\
\theta\left(  T,x\right)  =\dfrac{a}{x}.
\end{array}
\right.  \tag{5.7}%
\end{equation}

Once again, we know that the solution of $\left(  5.7\right)  $ will be of the form%

\begin{equation}
\theta\left(  t,x\right)  =\Pi\left(  t\right)  \dfrac{a}{x},\text{ for }%
t\in\left[  0,T\right]  , \tag{5.8}%
\end{equation}
where $\Pi\left(  \cdot\right)  \in C^{1}\left(  \left[  0,T\right]
,\mathbb{R}\right)  .$ By substituting in $\left(  5.7\right)  ,$ we get%
\begin{equation}
\left\{
\begin{array}
[c]{l}%
\Pi_{t}\left(  t\right)  +\dfrac{1}{a\lambda\left(  T-t\right)  }=0,\text{ for
}t\in\left[  0,T\right]  ,\\
\Pi\left(  T\right)  =1,
\end{array}
\right.  \tag{5.9}%
\end{equation}
which is explicitly solved by%
\[
\Pi\left(  t\right)  =1+\int_{t}^{T}\dfrac{1}{a\lambda\left(  T-r\right)
}dr,\text{ for }t\in\left[  0,T\right]  .
\]

In view of Theorem 4.1, the representation of the Nash equilibrium\ strategies
$\left(  4.3\right)  $-$\left(  4.4\right)  $ gives
\begin{align}
&  \hat{c}\left(  t\right)  =\left(  a\lambda\left(  T-t\right)  +\int_{t}%
^{T}\dfrac{\lambda\left(  T-t\right)  }{\lambda\left(  T-r\right)  }dr\right)
^{-1}\hat{X}\left(  t\right)  ,\text{ }dt-a.e.,\tag{5.10}\\
&  \hat{u}_{I}\left(  t\right)  =\Sigma\left(  t\right)  r\left(  t\right)
\hat{X}\left(  t\right)  ,\text{ }dt-a.e. \tag{5.11}%
\end{align}

This consumption--investment strategy determines a wealth process given by%
\[%
\begin{array}
[c]{ll}%
X\left(  t\right)  = & x_{0}+%
{\displaystyle\int_{0}^{t}}
\left\{  r_{0}\left(  s\right)  +r\left(  s\right)  ^{\top}\Sigma\left(
s\right)  r\left(  s\right)  -\left(  a\lambda\left(  T-s\right)  +\int
_{s}^{T}\dfrac{\lambda\left(  T-s\right)  }{\lambda\left(  T-r\right)
}dr\right)  ^{-1}\right\}  \hat{X}\left(  s\right)  ds\\
& \text{ \ \ \ \ }+%
{\displaystyle\int_{0}^{t}}
r\left(  s\right)  ^{\top}\Sigma\left(  s\right)  \sigma\left(  s\right)
\hat{X}\left(  s\right)  dW\left(  s\right)  .
\end{array}
\]

\subsection{Exponential utility function}

Next, we consider the case where $\varphi\left(  c\right)  =-\dfrac{e^{-\gamma
c}}{\gamma}$ and $h\left(  x\right)  =-a\dfrac{e^{-\gamma x}}{\gamma},$ with
$a,\gamma>0$. The terminal condition PDE $\left(  4.1\right)  $ becomes
\begin{equation}
\left\{
\begin{array}
[c]{l}%
\theta_{t}\left(  t,x\right)  +\theta_{x}\left(  t,x\right)  \left(
r_{0}\left(  t\right)  x-r\left(  t\right)  ^{\top}\Sigma\left(  t\right)
r\left(  t\right)  \dfrac{\theta\left(  t,x\right)  }{\theta_{x}\left(
t,x\right)  }-\dfrac{1}{\gamma}\ln\left(  \lambda\left(  T-t\right)
\theta\left(  t,x\right)  \right)  \right) \\
+\dfrac{1}{2}\theta_{xx}\left(  t,x\right)  r\left(  t\right)  ^{\top}%
\Sigma\left(  t\right)  r\left(  t\right)  \left(  \dfrac{\theta\left(
t,x\right)  }{\theta_{x}\left(  t,x\right)  }\right)  ^{2}+r_{0}\left(
t\right)  \theta\left(  t,x\right)  =0,\text{ }\left(  t,x\right)  \in\left[
0,T\right]  \times\mathbb{R} ,\\
\theta\left(  T,x\right)  =ae^{-\gamma x}.
\end{array}
\right.  \tag{5.12}%
\end{equation}

We try a solution of the form
\begin{equation}
\theta\left(  t,x\right)  =ae^{-\gamma\left(  \phi\left(  t\right)
x+\psi\left(  t\right)  \right)  },\text{ for }t\in\left[  0,T\right]  ,
\tag{5.13}%
\end{equation}
where $\phi\left(  \cdot\right)  $, $\psi\left(  \cdot\right)  \in
C^{1}\left(  \left[  0,T\right]  ,\mathbb{R}\right)  $ such that $\phi\left(
T\right)  =1$ and $\psi\left(  T\right)  =0.$ By substituting in $\left(
5.12\right)  $ we get%
\begin{align*}
\left\{  -\gamma\phi_{t}\left(  t\right)  +\gamma\phi\left(  t\right)
^{2}-\gamma\phi\left(  t\right)  r_{0}\left(  t\right)  \right\}  x-\dfrac
{1}{2}r\left(  t\right)  ^{\top}\Sigma\left(  t\right)  r\left(  t\right)   &
\\
-\gamma\psi_{t}\left(  t\right)  -\phi\left(  t\right)  \ln\left(
a\lambda\left(  T-t\right)  \right)  +\gamma\phi\left(  t\right)  \psi\left(
t\right)  +r_{0}\left(  t\right)  =0.  &
\end{align*}

This suggests that functions $\phi\left(  \cdot\right)  $ and $\psi\left(
\cdot\right)  $ should solve the following system of equations%
\begin{equation}
\left\{
\begin{array}
[c]{l}%
\phi_{t}\left(  t\right)  =-r_{0}\left(  t\right)  \phi\left(  t\right)
+\phi\left(  t\right)  ^{2},\text{ }t\in\left[  0,T\right]  ,\\
\psi_{t}\left(  t\right)  =-\dfrac{1}{\gamma}\phi\left(  t\right)  \ln\left(
a\lambda\left(  T-t\right)  \right)  +\phi\left(  t\right)  \psi\left(
t\right)  -\dfrac{1}{2\gamma}r\left(  t\right)  ^{\top}\Sigma\left(  t\right)
r\left(  t\right)  +\dfrac{1}{\gamma}r_{0}\left(  t\right)  ,\text{ }%
t\in\left[  0,T\right]  ,\\
\phi\left(  T\right)  =1,\text{ }\psi\left(  T\right)  =0,
\end{array}
\right.  \tag{5.14}%
\end{equation}
which is explicitly solvable for $t\in\left[  0,T\right]  ,$ by%
\begin{equation}
\phi\left(  t\right)  =\dfrac{e^{\int_{t}^{T}r_{0}\left(  \tau\right)  d\tau}%
}{1+\int_{t}^{T}e^{\int_{l}^{T}r_{0}\left(  \tau\right)  d\tau}dl}, \tag{5.15}%
\end{equation}
and%
\begin{equation}
\psi\left(  t\right)  =e^{-\int_{t}^{T}\phi\left(  \tau\right)  d\tau}\int
_{t}^{T}e^{\int_{l}^{T}\phi\left(  \tau\right)  d\tau}\left(  \dfrac{1}%
{\gamma}\phi\left(  l\right)  \ln\left(  \lambda\left(  T-l\right)  a\right)
+\dfrac{1}{2\gamma}r\left(  t\right)  ^{\top}\Sigma\left(  t\right)  r\left(
t\right)  -\dfrac{r_{0}\left(  l\right)  }{\gamma}\right)  dt. \tag{5.16}%
\end{equation}

The representation of the Nash equilibrium\ strategies $\left(  4.3\right)
$-$\left(  4.4\right)  $ gives%
\begin{align}
&  \hat{c}\left(  t\right)  =-\frac{1}{\gamma}\ln\left(  a\lambda\left(
T-t\right)  \right)  +\phi\left(  t\right)  \hat{X}\left(  t\right)
+\psi\left(  t\right)  ,\text{ }dt-a.e.\tag{5.17}\\
&  \hat{u}_{I}\left(  t\right)  =\frac{1}{\gamma}\Sigma\left(  t\right)
r\left(  t\right)  \phi\left(  t\right)  ^{-1},\text{ }dt-a.e. \tag{5.18}%
\end{align}

This consumption--investment strategy determines a wealth process given by%
\[%
\begin{array}
[c]{cl}%
X\left(  t\right)  = & x_{0}+%
{\displaystyle\int_{0}^{t}}
\left\{  \left(  r_{0}\left(  s\right)  -\phi\left(  s\right)  \right)
\hat{X}\left(  s\right)  +\dfrac{1}{\gamma}\left(  r\left(  s\right)  ^{\top
}\Sigma\left(  s\right)  r\left(  s\right)  \phi\left(  s\right)  ^{-1}%
-\ln\left(  a\lambda\left(  T-s\right)  \right)  \right)  \right. \\
& \text{ \ \ \ \ \ \ \ \ \ \ \ \ \ \ \ }\left.  -\psi\underset{}{\left(
s\right)  }\right\}  ds+%
{\displaystyle\int_{0}^{t}}
\dfrac{1}{\gamma}\phi\left(  s\right)  ^{-1}r\left(  s\right)  ^{\top}%
\Sigma\left(  s\right)  \sigma\left(  s\right)  dW\left(  s\right)  ,\text{
}t\in\left[  0,T\right]  ,
\end{array}
\]

The above solution is comparable with the ones obtained in Mar\'{\i}n-Solano
\& Navas \cite{Solano} by solving an extended Hamilton--Jacobi--Bellman (HJB) equations.

\section{Special discount function}

As well documented in \cite{Solano}, an agent making a decision at time $t$ is
usually called the $t$-agent, and can act in two different ways: naive and
sophisticated. Naive agents take decisions without taking into account that
their preferences will change in the near future, and then any $t$-agent will
solve the problem as a standard optimal control problem with initial condition
$X(t)=x_{t}$ and his decision will be in general time-inconsistent. In order
to obtain a time consistent strategy, the $t$-agent should be sophisticated,
in the sense of taking into account the preferences of all the $s$-agents, for
$s\in\left[  t,T\right]  $. Therefore, the approach to handle the time
inconsistency in dynamic decision making problems is by considering
time-inconsistent problems as non-cooperative games with a continuous number
of players, in which decisions at every instant of time are selected. The
solution to the problem of the agent with non-constant discounting should be
constructed by looking for the sub-game perfect equilibria of the associated
game with an infinite number of $t$-agents. In \cite{Solano} the authors
looked for a solution of a sophisticated agent to the modified HJB (which is
not a partial differential equation due to the presence of a non-local term).
Then, they need to define the Markov equilibrium strategies, while in our
work, and different from \cite{Solano}, we use the open-loop equilibrium
strategies. This is a significant difference which leads to obtain an
important change in the results.

\subsection{Exponential discounting with constant discount rate (classical
model)}

At first, we consider the standard exponential discount function
$\lambda\left(  t\right)  =e^{-\delta_{0}t}$, $t\in\left[  0,T\right]  $,
where $\delta_{0}>0$ is a constant representing the discount rate. In this
case, our equilibrium solution for the three cases become \vspace{0.5cm}

1) \textbf{Logarithmic utility}%

\begin{align*}
&  \hat{c}(t)=\frac{1}{ae^{-\left(  T-t\right)  \delta_{0}} +\int_{s}%
^{T}e^{-\left(  l-t\right)  \delta_{0}}dl}\hat{X}\left(  t\right)  ,\text{
}dt-a.e.,\\
&  \hat{u}_{I}(t)=\Sigma\left(  t\right)  r\left(  t\right)  \hat{X}\left(
t\right)  ,\text{ }dt-a.e.
\end{align*}

2) \textbf{Power utility}%
\begin{align*}
&  \hat{c}\left(  t\right)  \mathbf{=}\left(  ae^{-\left(  T-t\right)
\delta_{0}}\right)  ^{\frac{1}{\gamma-1}}\frac{e^{\int_{t}^{T}\dfrac{K\left(
\tau\right)  }{\gamma-1}d\tau}}{\left(  1+\int_{t}^{T}\left(  ae^{-\left(
T-\tau\right)  \delta_{0}}\right)  ^{\frac{1}{\gamma-1}}e^{\int_{\tau}%
^{T}\dfrac{K\left(  l\right)  }{\gamma-1}dl}d\tau\right)  }\hat{X}\left(
t\right)  ,\text{ }dt-a.e.,\\
&  \hat{u}_{I}(t)=\Sigma\left(  t\right)  r\left(  t\right)  \frac{\hat
{X}\left(  t\right)  }{\left(  1-\gamma\right)  },\text{ }dt-a.e.
\end{align*}

3) \textbf{Exponential utility}%
\begin{align*}
&  \hat{c}\left(  t\right)  \mathbf{=}-\frac{1}{\gamma}\ln\left(  ae^{-\left(
T-t\right)  \delta_{0}}\right)  +\phi\left(  t\right)  \hat{X}\left(
t\right)  +\psi\left(  t\right)  ,\text{ }dt-a.e.,\\
&  \hat{u}_{I}\left(  t\right)  =\Sigma\left(  t\right)  r\left(  t\right)
\frac{1}{\gamma\phi\left(  t\right)  },\text{ }dt-a.e.
\end{align*}
where $K\left(  \cdot\right)  ,$ $\phi\left(  \cdot\right)  $ are given by
$\left(  5.2\right)  $ and $\left(  5.15\right)  ,$ respectively, and
\[
\psi\left(  t\right)  =\dfrac{1}{\gamma}e^{-\int_{t}^{T}\phi\left(
\tau\right)  d\tau} \int_{t}^{T}e^{\int_{l}^{T}\phi\left(  \tau\right)  d\tau
}\left(  \phi\left(  l\right)  \ln\left(  e^{-\left(  T-l\right)  \delta_{0}%
}a\right)  +\dfrac{1}{2}r\left(  l\right)  ^{\top}\Sigma\left(  l\right)
r\left(  l\right)  -r_{0}\left(  l\right)  \right)  dl.
\]

Notice that our solutions given above coincide with the optimal solutions of
classical Merton portfolio problem (see e.g.\cite{Solano} in the case with
constant discount rate). This confirms the well-known fact that the
time-consistent equilibrium strategy for an exponential discount function is
nothing but the optimal strategy. A relevant remark is that the portfolio rule
is independent of the discount factor, and it is the same for a
non-exponential discount function.

\subsection{Exponential discounting with non constant discount rate\newline%
(Karp's model)}

Now, following Karp \cite{Karp}, let us assume that the instantaneous discount
rate is non-constant, but a continuous and positive function of time
$\delta\left(  l\right)  $, for $l\in\left[  0,T\right]  .$ Impatient agents
will be characterized by a non-increasing discount rate $\delta\left(
\cdot\right)  $. The discount factor used to evaluate a payoff at times $\tau$
$\geq$ 0, is given by%

\begin{equation}
\lambda\left(  \tau\right)  =e^{-\int_{0}^{\tau}\delta\left(  l\right)  dl}.
\tag{6.1}%
\end{equation}

In this case, the objective is exactly the same as Mar\'{\i}n-Solano and Navas
\cite{Solano}, in which the equilibrium is however defined within the class of
feedback controls. In \cite{Solano}, the (feedback) equilibrium
consumption-investment solutions (also called the sophisticated
consumption-investment strategies) are summarized as \vspace{0.5cm}

1) \textbf{Logarithmic utility}$\ $%
\begin{align}
&  \hat{c}\left(  t\right)  =\frac{1}{ae^{-\int_{0}^{T-t}\delta\left(
\tau\right)  d\tau}+\int_{t}^{T}e^{-\int_{0}^{s-l}\delta\left(  \tau\right)
d\tau}dl}\hat{X}\left(  t\right)  ,\text{ }dt-a.e.,\tag{6.2}\\
&  \hat{u}_{I}(t)=\Sigma\left(  t\right)  r\left(  t\right)  \hat{X}\left(
t\right)  ,\text{ }dt-a.e. \tag{6.3}%
\end{align}

2) \textbf{Power utility}%
\begin{align}
&  \hat{c}\left(  t\right)  =\left(  \alpha\left(  t\right)  \right)
^{\frac{1}{\gamma-1}}\hat{X}\left(  t\right)  ,\text{ }dt-a.e.,\tag{6.4}\\
&  \hat{u}_{I}(t)=\Sigma\left(  t\right)  r\left(  t\right)  \frac{\hat
{X}\left(  t\right)  }{\left(  1-\gamma\right)  },\text{ }dt-a.e. \tag{6.5}%
\end{align}
where $\alpha\left(  \cdot\right)  $ is the solution of the
integro-differential equation,%

\begin{equation}
\left\{
\begin{array}
[c]{l}%
\alpha_{t}\left(  t\right)  -\left(  \delta\left(  T-t\right)  -K\left(
t\right)  \right)  \alpha\left(  t\right)  +\left(  1-\gamma\right)
\alpha\left(  t\right)  ^{\frac{\gamma}{1-\gamma}}\\
-{\displaystyle\int_{t}^{T}}e^{-\int_{0}^{s-t}\delta\left(  l\right)
dl}\left(  \delta\left(  s-t\right)  -\delta\left(  T-t\right)  \right)
\alpha\left(  s\right)  ^{\frac{\gamma}{1-\gamma}}e^{\gamma\int_{t}^{s}%
\Delta\left(  \tau\right)  d\tau}ds=0,\\
\alpha\left(  T\right)  =a.
\end{array}
\right.  \tag{6.6}%
\end{equation}
with $K\left(  t\right)  $ given by $\left(  5.2\right)  $ and%
\[
\Delta\left(  \tau\right)  =r_{0}\left(  \tau\right)  +\dfrac{1}{\left(
1-\gamma\right)  }r\left(  \tau\right)  ^{\top}\Sigma\left(  \tau\right)
r\left(  \tau\right)  -\alpha\left(  \tau\right)  ^{\frac{1}{1-\gamma}}.
\]

3) \textbf{Exponential utility}%
\begin{align}
&  \hat{c}\left(  t\right)  =\phi\left(  t\right)  \hat{X}\left(  t\right)
+C\left(  t\right)  -\frac{\ln\left(  \gamma a\phi\left(  t\right)  \right)
}{\gamma},\text{ }dt-a.e.,\tag{6.7}\\
&  \hat{u}_{I}(t)=\Sigma\left(  t\right)  r\left(  t\right)  \frac{1}%
{\gamma\phi\left(  t\right)  },\text{ }dt-a.e. \tag{6.8}%
\end{align}
where $\phi\left(  \cdot\right)  $ is given by $\left(  5.15\right)  $ and
$C\left(  \cdot\right)  $ satisfies the following very complicated
integro-differential equation,%

\begin{equation}
\left\{
\begin{array}
[c]{l}%
C_{t}\left(  t\right)  -C\left(  t\right)  \phi\left(  t\right)  +\dfrac
{1}{\gamma}\phi\left(  t\right)  \ln\left(  a\gamma\phi\left(  t\right)
\right)  +\dfrac{1}{2\gamma}r\left(  t\right)  ^{\top}\Sigma\left(  t\right)
r\left(  t\right) \\
+\dfrac{1}{\gamma}\left\{  \delta\left(  T-t\right)  -\phi\left(  t\right)
-\mathcal{K}\left(  C\left(  t\right)  ,t\right)  \right\}  =0,\\
C\left(  T\right)  =0,
\end{array}
\right.  \tag{6.9}%
\end{equation}
where
\begin{align}
&  \mathcal{K}\left(  C\left(  t\right)  ,t\right)  =-\mathbb{E} \left[
\int_{t}^{T}e^{-\int_{0}^{s-t}\delta\left(  l\right)  dl}\left\{
\delta\left(  s-t\right)  -\delta\left(  T-t\right)  \right\}  \phi\left(
t\right)  \right. \nonumber\\
&  \text{
\ \ \ \ \ \ \ \ \ \ \ \ \ \ \ \ \ \ \ \ \ \ \ \ \ \ \ \ \ \ \ \ \ \ \ \ \ }%
\times\left.  e^{-\gamma\left\{  C\left(  s\right)  -C\left(  t\right)
+\int_{t}^{s}\phi\left(  \tau\right)  Z\left(  \tau\right)  d\tau+\int_{t}%
^{s}\frac{1}{\gamma}r\left(  \tau\right)  ^{\top}\Sigma\left(  \tau\right)
\mathbf{\sigma}\left(  \tau\right)  dW\left(  \tau\right)  \right\}
}ds\right]  , \tag{6.10}%
\end{align}
with
\[
Z\left(  \tau\right)  =\dfrac{1}{\gamma\phi\left(  \tau\right)  }r\left(
\tau\right)  ^{\top}\Sigma\left(  \tau\right)  r\left(  \tau\right)  -C\left(
\tau\right)  +\frac{1}{\gamma}\ln\left(  \gamma a\phi\left(  \tau\right)
\right)  .
\]

Our (open-loop) equilibrium solutions reduce to \vspace{0.5cm}

1) \textbf{Logarithmic utility}$\ $%
\begin{align}
&  \hat{c}\left(  t\right)  =\frac{1}{ae^{-\int_{0}^{T-t}\delta\left(
\tau\right)  d\tau}+\int_{t}^{T}e^{-\int_{T-l}^{T-s}\delta\left(  \tau\right)
d\tau}dl}\hat{X}\left(  t\right)  ,\text{ }dt-a.e.,\tag{6.11}\\
&  \hat{u}_{I}(t)=\Sigma\left(  t\right)  r\left(  t\right)  \hat{X}\left(
t\right)  ,\text{ }dt-a.e. \tag{6.12}%
\end{align}

2) \textbf{Power utility}%
\begin{align}
&  \hat{c}\left(  t\right)  =\frac{\left(  ae^{-\int_{0}^{T-t}\delta\left(
\tau\right)  d\tau}\right)  ^{\frac{1}{\gamma-1}}e^{\int_{s}^{T}%
\dfrac{K\left(  \tau\right)  }{\gamma-1}d\tau}}{\left(  1+\int_{t}^{T}\left(
ae^{-\int_{0}^{T-\tau}\delta\left(  \tau\right)  d\tau}\right)  ^{\frac
{1}{\gamma-1}}e^{\int_{\tau}^{T}\dfrac{K\left(  l\right)  }{\gamma-1}dl}%
d\tau\right)  }\hat{X}\left(  t\right)  ,\text{ }dt-a.e.,\tag{6.13}\\
&  \hat{u}_{I}(t)=\Sigma\left(  t\right)  r\left(  t\right)  \frac{\hat
{X}\left(  t\right)  }{\left(  1-\gamma\right)  },\text{ }dt-a.e. \tag{6.14}%
\end{align}

3) \textbf{Exponential utility}%

\begin{align}
&  \hat{c}\left(  t\right)  =-\frac{1}{\gamma}\ln\left(  ae^{-\int_{0}%
^{T-t}\delta\left(  \tau\right)  d\tau}\right)  +\phi\left(  t\right)  \hat
{X}\left(  t\right)  +\psi\left(  t\right)  ,\text{ }dt-a.e.,\tag{6.15}\\
&  \hat{u}_{I}(t)=\Sigma\left(  t\right)  r\left(  t\right)  \frac{1}%
{\gamma\hat{X}\left(  t\right)  \phi\left(  t\right)  },\text{ }dt-a.e.
\tag{6.16}%
\end{align}
where $K\left(  \cdot\right)  ,$ $\phi\left(  \cdot\right)  $ are given by
$\left(  5.2\right)  $ and $\left(  5.15\right)  ,$ respectively, and%

\[
\psi\left(  t\right)  =e^{-\int_{t}^{T}\phi\left(  \tau\right)  d\tau}\int
_{t}^{T}e^{\int_{l}^{T}\phi\left(  \tau\right)  d\tau}\left(  \dfrac{1}%
{\gamma}\phi\left(  l\right)  \ln\left(  e^{-\int_{0}^{T-t}\delta\left(
\tau\right)  d\tau}a\right)  +\dfrac{1}{2\gamma}r\left(  l\right)  ^{\top
}\Sigma\left(  l\right)  r\left(  l\right)  -\dfrac{r_{0}\left(  l\right)
}{\gamma}\right)  dl.
\]

\begin{remark}
Comparing the results of this special case with our solutions, we find the
following facts: The equilibrium proportion investment strategies coincide in
the three cases. The consumption strategies are different in the three cases.
Moreover, our equilibrium consumption strategies are well defined and
explicitly given, while in \cite{Solano}, equilibrium consumption strategies
in the case of Power utility as well as in the case of Exponential utility,
are obtained via a very complicated integro-differential equations, whose
unique solvability are not established.
\end{remark}

\section{Appendix}

Following \cite{Huetal2011}, we derive the proof of Proposition 3.3 by means
of the duality analysis. Moreover, since our objective function is not in
quadratic form, we need to adapt the results obtained in \cite{Huetal2011}
according to our control problem which concerns a general and non necessary
quadratic utility maximization. \vspace{0.5cm}

\textbf{Proof of Proposition 3.3.} The estimates $\left(  3.15\right)
-(3.17)$ follow from Theorem 4.4 in \cite{YongZhou}. Moreover the following
expansion holds for the objective functional%
\begin{align}
&  J\left(  t,\hat{X}\left(  t\right)  ,u^{\varepsilon}\left(  \cdot\right)
\right)  -J\left(  t,\hat{X}\left(  t\right)  ,\hat{u}\left(  \cdot\right)
\right) \nonumber\\
&  =\mathbb{E}^{t}\left[  {\displaystyle\int_{t}^{T}} \lambda\left(
s-t\right)  \left(  \varphi\left(  \Gamma^{\top}u^{\varepsilon}\left(
s\right)  \right)  -\varphi\left(  \Gamma^{\top}\hat{u}\left(  s\right)
\right)  \right)  ds+\lambda\left(  T-t\right)  \left(  h\left(
X^{\varepsilon}\left(  T\right)  \right)  -h\left(  \hat{X}\left(  T\right)
\right)  \right)  \right]  . \tag{A.1.1}%
\end{align}

Now, applying the second order Taylor-Lagrange expansion to $\varphi\left(
\Gamma^{\top}u^{\varepsilon}\left(  s\right)  \right)  -\varphi\left(
\Gamma^{\top}u\left(  s\right)  \right)  ,$ we find%
\begin{align}
&  \varphi\left(  \Gamma^{\top}u^{\varepsilon}\left(  s\right)  \right)
-\varphi\left(  \Gamma^{\top}\hat{u}\left(  s\right)  \right) \nonumber\\
&  =\varphi\left(  \Gamma^{\top}\left(  \hat{u}\left(  s\right)  +v1_{\left[
t,t+\varepsilon\right)  }\right)  \right)  -\varphi\left(  \Gamma^{\top}%
\hat{u}\left(  s\right)  \right)  ,\nonumber\\
&  =\left\{  \left\langle \varphi_{c}\left(  \Gamma^{\top}\hat{u}\left(
s\right)  \right)  \Gamma,v\right\rangle +\frac{1}{2}\left\langle \varphi
_{cc}\left(  \Gamma^{\top}\left(  \hat{u}\left(  s\right)  +\theta v1_{\left[
t,t+\varepsilon\right)  }\right)  \right)  \Gamma\Gamma^{\top}v,v\right\rangle
\right\}  1_{\left[  t,t+\varepsilon\right)  }.\nonumber\\
&  =\left\langle \varphi_{c}\left(  \Gamma^{\top}\hat{u}\left(  s\right)
\right)  \Gamma,v\right\rangle 1_{\left[  t,t+\varepsilon\right)  }+\frac
{1}{2}\left\langle \varphi_{cc}\left(  \Gamma^{\top}\hat{u}\left(  s\right)
\right)  \Gamma\Gamma^{\top}v,v\right\rangle 1_{\left[  t,t+\varepsilon
\right)  }\nonumber\\
&  +\frac{1}{2}\left\{  \left\langle \left(  \varphi_{cc}\left(  \Gamma^{\top
}\left(  \hat{u}\left(  s\right)  +\theta v1_{\left[  t,t+\varepsilon\right)
}\right)  \right)  -\varphi_{cc}\left(  \Gamma^{\top}\hat{u}\left(  s\right)
\right)  \right)  \Gamma\Gamma^{\top}v,v\right\rangle \right\}  1_{\left[
t,t+\varepsilon\right)  }. \tag{A.1.2}%
\end{align}

Notice that%
\begin{align*}
&  \mathbb{E}^{t}\left[  {\displaystyle\int_{t}^{T}} \lambda\left(
s-t\right)  \left\{  \left\langle \left(  \varphi_{cc}\left(  \Gamma^{\top
}\left(  \hat{u}\left(  s\right)  +\theta v1_{\left[  t,t+\varepsilon\right)
}\right)  \right)  -\varphi_{cc}\left(  \Gamma^{\top}\hat{u}\left(  s\right)
\right)  \right)  \Gamma\Gamma^{\top}v,v\right\rangle \right\}  1_{\left[
t,t+\varepsilon\right)  }ds\right] \\
&  \leq C\left\Vert v\right\Vert ^{2}\mathbb{E}^{t}\left[  \left(
{\displaystyle\int_{t}^{T}} \theta\left\vert v1_{\left[  t,t+\varepsilon
\right)  }\right\vert ds\right)  ^{2}\right]  ^{\frac{1}{2}}\mathbb{E}%
^{t}\left[  \left(  {\displaystyle\int_{t}^{T}} 1_{\left[  t,t+\varepsilon
\right)  }ds\right)  ^{2}\right]  ^{\frac{1}{2}}\\
&  =C\theta\left\Vert v\right\Vert ^{3}\varepsilon^{2}=o\left(  \varepsilon
\right)  ,
\end{align*}
where $\left\Vert v\right\Vert $ denote the Euclidean norm of $v$, that is
$\left\Vert v\right\Vert ^{2}=:\left\langle v,v\right\rangle .$ Then we obtain%
\begin{align*}
&  \mathbb{E}^{t}\left[  {\displaystyle\int_{t}^{T}} \lambda\left(
s-t\right)  \left(  \varphi\left(  \Gamma^{\top}u^{\varepsilon}\left(
s\right)  \right)  -\varphi\left(  \Gamma^{\top}\hat{u}\left(  s\right)
\right)  \right)  ds\right] \\
&  =\mathbb{E}^{t}\left[  {\displaystyle\int_{t}^{T}} \lambda\left(
s-t\right)  \left\{  \left\langle \varphi_{c}\left(  \Gamma^{\top}\hat
{u}\left(  s\right)  \right)  \Gamma,v\right\rangle +\frac{1}{2}\left\langle
\varphi_{cc}\left(  \Gamma^{\top}\hat{u}\left(  s\right)  \right)
\Gamma\Gamma^{\top}v,v\right\rangle \right\}  1_{\left[  t,t+\varepsilon
\right)  }\right] \\
&  +o\left(  \varepsilon\right)  .
\end{align*}

Noting that, by the second order Taylor-Lagrange expansion, see e.g.
\cite{YongZhou}, we have for some constant $L>0,$
\begin{align*}
&  \mathbb{E}^{t}\left[  h\left(  \hat{X}\left(  T\right)  +y^{\varepsilon
,v}\left(  T\right)  +z^{\varepsilon,v}\left(  T\right)  \right)  -h\left(
\hat{X}\left(  T\right)  \right)  \right] \\
&  =\mathbb{E}^{t}\left[  h_{x}\left(  \hat{X}\left(  T\right)  \right)
\left(  y^{\varepsilon,v}\left(  T\right)  +z^{\varepsilon,v}\left(  T\right)
\right)  +\frac{1}{2}h_{xx}\left(  \hat{X}\left(  T\right)  \right)  \left(
y^{\varepsilon,v}\left(  T\right)  +z^{\varepsilon,v}\left(  T\right)
\right)  ^{2}\right]  +o\left(  \varepsilon\right)  ,
\end{align*}
then the following expansion holds for the objective functional%
\begin{align}
&  J\left(  t,\hat{X}\left(  t\right)  ,u^{\varepsilon}\left(  \cdot\right)
\right)  -J\left(  t,\hat{X}\left(  t\right)  ,\hat{u}\left(  \cdot\right)
\right) \tag{A.1.3}\\
&  =\mathbb{E}^{t}\left[  {\displaystyle\int_{t}^{T}} \lambda\left(
s-t\right)  \left\{  \left\langle \varphi_{c}\left(  \Gamma^{\top}\hat
{u}\left(  s\right)  \right)  \Gamma,v\right\rangle +\frac{1}{2}\left\langle
\varphi_{cc}\left(  \Gamma^{\top}\hat{u}\left(  s\right)  \right)
\Gamma\Gamma^{\top}v,v\right\rangle \right\}  1_{\left[  t,t+\varepsilon
\right)  }ds\right. \nonumber\\
&  +\left.  \lambda\left(  T-t\right)  \left(  h_{x}\left(  \hat{X}\left(
T\right)  \right)  \left(  y^{\varepsilon,v}\left(  s\right)  +z^{\varepsilon
,v}\left(  s\right)  \right)  +\frac{1}{2}h_{xx}\left(  \hat{X}\left(
T\right)  \right)  \left(  y^{\varepsilon,v}\left(  s\right)  +z^{\varepsilon
,v}\left(  s\right)  \right)  ^{2}\right)  \right]  +o\left(  \varepsilon
\right)  .\nonumber
\end{align}

Notice that%
\begin{align*}
&  \lambda\left(  T-t\right)  \left(  h_{x}\left(  \hat{X}\left(  T\right)
\right)  \left(  y^{\varepsilon,v}\left(  s\right)  +z^{\varepsilon,v}\left(
s\right)  \right)  +\frac{1}{2}h_{xx}\left(  \hat{X}\left(  T\right)  \right)
\left(  y^{\varepsilon,v}\left(  s\right)  +z^{\varepsilon,v}\left(  s\right)
\right)  ^{2}\right) \\
&  =p\left(  T;t\right)  \left(  y^{\varepsilon,v}\left(  s\right)
+z^{\varepsilon,v}\left(  s\right)  \right)  +\frac{1}{2}P\left(  s;t\right)
\left(  y^{\varepsilon,v}\left(  s\right)  +z^{\varepsilon,v}\left(  s\right)
\right)  ^{2}.
\end{align*}

Now, by applying Ito's formula to $s\mapsto p\left(  s;t\right)  \left(
y^{\varepsilon,v}\left(  s\right)  +z^{\varepsilon,v}\left(  s\right)
\right)  $ on $\left[  t,T\right]  $, we get%
\begin{equation}
\mathbb{E}^{t}\left[  p\left(  T;t\right)  \left(  y^{\varepsilon,v}\left(
T\right)  +z^{\varepsilon,v}\left(  T\right)  \right)  \right]  =\mathbb{E}%
^{t}\left[  {\displaystyle\int_{t}^{t+\varepsilon}} \left\{  v^{\top}B\left(
s\right)  p\left(  s;t\right)  +v^{\top}D\left(  s\right)  \widetilde
{q}\left(  s;t\right)  \right\}  ds\right]  . \tag{A.1.4}%
\end{equation}

Again, by applying Ito's formula to $s\mapsto P\left(  s;t\right)  \left(
y^{\varepsilon,v}\left(  s\right)  +z^{\varepsilon,v}\left(  s\right)
\right)  ^{2}$ on $\left[  t,T\right]  ,$ we get%
\begin{equation}%
\begin{array}
[c]{l}%
\mathbb{E}^{t}\left[  P\left(  T;t\right)  \left(  y^{\varepsilon,v}\left(
T\right)  +z^{\varepsilon,v}\left(  T\right)  \right)  ^{2}\right] \\
=\mathbb{E}^{t}\left[  {\displaystyle\int_{t}^{t+\varepsilon}} \left\{
2v^{\top}\left(  y^{\varepsilon,v}\left(  s\right)  +z^{\varepsilon,v}\left(
s\right)  \right)  \left(  B\left(  s\right)  P\left(  s,t\right)  +D\left(
s\right)  \widetilde{Q}\left(  s,t\right)  \right)  \right.  \right. \\
\text{ \ \ \ \ \ \ \ \ \ \ \ \ \ \ \ \ \ \ \ \ }\left.  \left.  +v^{\top
}\left(  D\left(  s\right)  D\left(  s\right)  ^{\top}\right)  vP\left(
s,t\right)  \right\}  ds\right]  ,
\end{array}
\tag{A.1.5}%
\end{equation}
where $\widetilde{Q}\left(  s;t\right)  =\left(  0,Q\left(  s;t\right)
^{\top}\right)  ^{\top}.$ On the other hand, we conclude from \textbf{(H1)}
together with $\left(  3.17\right)  $ that%
\begin{equation}
\mathbb{E}^{t}\left[  {\displaystyle\int_{t}^{t+\varepsilon}} \left(
y^{\varepsilon,v}\left(  s\right)  +z^{\varepsilon,v}\left(  s\right)
\right)  \left(  B\left(  s\right)  P\left(  s,t\right)  +D\left(  s\right)
\widetilde{Q}\left(  s,t\right)  \right)  ds\right]  =o\left(  \varepsilon
\right)  . \tag{A.1.6}%
\end{equation}

By taking $\left(  A.1.4\right)  ,\left(  A.1.5\right)  $ and $\left(
A.1.6\right)  $ in $\left(  A.1.3\right)  ,$ it follows that%
\begin{align*}
&  J\left(  t,\hat{X}\left(  t\right)  ,u^{\varepsilon}\left(  \cdot\right)
\right)  -J\left(  t,\hat{X}\left(  t\right)  ,\hat{u}\left(  \cdot\right)
\right) \\
&  =\mathbb{E}^{t}\left[  {\displaystyle\int_{t}^{t+\varepsilon}} \left\{
\left\langle B\left(  s\right)  p\left(  s;t\right)  +D\left(  s\right)
\widetilde{q}\left(  s;t\right)  +\lambda\left(  s-t\right)  \varphi
_{c}\left(  \Gamma^{\top}\hat{u}\left(  s\right)  \right)  \mathbf{1}_{\left[
t,t+\varepsilon\right)  }\Gamma,v\right\rangle \right.  \right. \\
&  \text{ \ \ \ \ \ \ \ \ \ \ \ \ }\left.  +\left.  \frac{1}{2}\left\langle
\left(  \varphi_{cc}\left(  \left\langle \Gamma\mathbf{,}\hat{u}\left(
s\right)  \right\rangle \right)  \Gamma\Gamma^{\top}+P\left(  s,t\right)
D\left(  s\right)  D\left(  s\right)  ^{\top}\right)  v,v\right\rangle
\right\}  ds\right]  +o\left(  \varepsilon\right)  ,
\end{align*}
which is equivalent to $(3.18).$\eop

\textbf{Proof of Lemma 3.4. } We set up%

\[
\alpha\left(  s\right)  =e^{-\int_{s}^{T} r_{0}\left(  \tau\right)  d\tau}.
\]
Now we define, for $t\in\left[  0,T\right]  $ and $s\in\left[  t,T\right]  ,$%

\[
\left(  \bar{p}\left(  s;t\right)  ,\bar{q}\left(  s;t\right)  \right)
:=\frac{\alpha\left(  s\right)  }{\lambda\left(  T-t\right)  }\left(  p\left(
s;t\right)  ,q\left(  s;t\right)  \right)  .
\]

Then, for any $t\in\left[  0,T\right]  ,$ in the interval $\left[  t,T\right]
,$ the pair $\left(  \bar{p}\left(  \cdot;t\right)  ,\bar{q}\left(
\cdot;t\right)  \right)  $ satisfies%
\begin{equation}
\left\{
\begin{array}
[c]{l}%
d\bar{p}\left(  s;t\right)  =\bar{q}\left(  s;t\right)  ^{\top}dW\left(
s\right)  ,\text{\ }s\in\left[  t,T\right]  ,\\
\bar{p}\left(  T;t\right)  =h_{x}\left(  \hat{X}\left(  T\right)  \right)  ,
\end{array}
\right.  \tag{A.2.1}%
\end{equation}
Moreover, it is clear that from the uniqueness of solutions to $\left(
A.2.1\right)  $, we have the equality $\left(  \bar{p}\left(  s;t_{1}\right)
,\bar{q}\left(  s;t_{1}\right)  \right)  =\left(  \bar{p}\left(
s;t_{2}\right)  ,\bar{q}\left(  s;t_{2}\right)  \right)  ,$ for any
$t_{1},t_{2},s\in\left[  0,T\right]  $ such that $0<t_{1}<t_{2}<s<T.$ Hence,
the solution $\left(  \bar{p}\left(  \cdot;t\right)  ,\bar{q}\left(
\cdot;t\right)  \right)  $ does not depend on the variable $t$ and this allows
us to denote the solution of $\left(  A.2.1\right)  $ by $\left(  \bar
{p}\left(  \cdot\right)  ,\bar{q}\left(  \cdot\right)  \right)  .$

We have then, for any $t\in\left[  0,T\right]  ,$ and $s\in\left[  t,T\right]
,$%
\begin{equation}
\left(  p\left(  s;t\right)  ,q\left(  s;t\right)  \right)  =\lambda\left(
T-t\right)  \alpha\left(  s\right)  ^{-1}\left(  \bar{p}\left(  s\right)
,\bar{q}\left(  s\right)  \right)  . \tag{A.2.2}%
\end{equation}

Now using $\left(  A.2.2\right)  $ we have, under (\textbf{H2)}, for any
$t\in\left[  0,T\right]  $ and $s\in\left[  t,T\right]  ,$
\begin{equation}%
\begin{array}
[c]{l}%
\left\vert p\left(  s;t\right)  -p\left(  s;s\right)  \right\vert =\left\vert
\left(  \lambda\left(  T-t\right)  -\lambda\left(  T-s\right)  \right)
\alpha\left(  s\right)  ^{-1}\bar{p}\left(  s\right)  \right\vert \\
\text{ \ \ \ \ \ \ \ \ \ \ \ \ \ \ \ \ \ \ \ \ \ \ \ } \leq C \left\vert
s-t\right\vert \left\vert \alpha\left(  s\right)  ^{-1}\bar{p}\left(
s\right)  \right\vert ,
\end{array}
\tag{A.2.3}%
\end{equation}
and%

\begin{equation}
\left\vert q\left(  s;t\right)  -q\left(  s;s\right)  \right\vert \leq C
\left\vert s-t\right\vert \left\vert \alpha\left(  s\right)  ^{-1}\bar
{q}\left(  s\right)  \right\vert . \tag{A.2.4}%
\end{equation}

From which, we have for any $t\in\left[  0,T\right]  ,$%

\begin{align*}
&  \underset{\varepsilon\downarrow0}{\lim}\dfrac{1}{\varepsilon}\mathbb{E}%
^{t}\left[  {\displaystyle\int_{t}^{t+\varepsilon}} \left\vert \mathcal{H}%
\left(  s;t\right)  -\mathcal{H}\left(  s;s\right)  \right\vert ds\right] \\
&  \leq C \underset{\varepsilon\downarrow0}{\lim}\dfrac{1}{\varepsilon
}\mathbb{E}^{t}\left[  {\displaystyle\int_{t}^{t+\varepsilon}} \left\vert
s-t\right\vert \left\vert \bar{p}\left(  s\right)  +\bar{q}\left(  s\right)
+\Gamma^{\top}\varphi_{c}\left(  \Gamma^{\top}\hat{u}\left(  s\right)
\right)  \right\vert ds\right]  ,\\
&  \leq C \underset{\varepsilon\downarrow0}{\lim}\mathbb{E}^{t}\left[
{\displaystyle\int_{t}^{t+\varepsilon}} \left\vert \bar{p}\left(  s\right)
+\bar{q}\left(  s\right)  +\Gamma^{\top}\varphi_{c}\left(  \Gamma^{\top}%
\hat{u}\left(  s\right)  \right)  \right\vert ds\right]  ,\\
&  =0.
\end{align*}

Thus
\begin{equation}
\underset{\varepsilon\downarrow0}{\lim}\dfrac{1}{\varepsilon}\mathbb{E}%
^{t}\left[  {\displaystyle\int_{t}^{t+\varepsilon}} \mathcal{H}\left(
s;t\right)  ds\right]  =\lim_{\varepsilon\downarrow0}\dfrac{1}{\varepsilon
}\mathbb{E}^{t}\left[  {\displaystyle\int_{t}^{t+\varepsilon}} \mathcal{H}%
\left(  s;s\right)  ds\right]  . \tag{A.2.5}%
\end{equation}

From the above equality, it is clear that if (ii) holds, then,%

\[
\underset{\varepsilon\downarrow0}{\lim}\dfrac{1}{\varepsilon} \mathbb{E}%
^{t}\left[  {\displaystyle\int_{t}^{t+\varepsilon}}\mathcal{H}\left(
s;t\right)  ds\right]  =0.\text{ }d\mathbb{P}-a.s,
\]

Conversely, according to Lemma 3.4 in \cite{Huetal2015}, if (i) holds, then,%

\[
\mathcal{H}\left(  s;s\right)  =0,\text{ }d\mathbb{P}-a.s,\text{ }ds-a.e.
\]

This completes the proof.\eop

\end{document}